\documentclass[12pt]{amsart}

\usepackage[latin1]{inputenc}
\usepackage{subfigure}
\usepackage{amsmath}
\usepackage{color}
\usepackage{amsmath, amsfonts, amssymb, hyperref, enumerate}
\usepackage{graphics}

\newcommand{\N}{\mathbb{N}}

\newcommand{\X}{\mathcal{X}}
\definecolor{gr75}{gray}{0.75}

\numberwithin{equation}{section}
\newcommand{\comment}[1]{\vspace{5 mm}\par \noindent
\marginpar{\textsc{Comment}}
\framebox{\begin{minipage}[c]{0.95 \textwidth}
 #1 \end{minipage}}\vspace{5 mm}\par}

\renewcommand{\comment}[1]{}

\newcommand{\ms}{\begin{math}}
\newcommand{\me}{\end{math}}

\newlength\cellsize 
\savebox2{%
\begin{picture}(15,15)
\put(0,0){\line(1,0){15}}
\put(0,0){\line(0,1){15}}
\put(15,0){\line(0,1){15}}
\put(0,15){\line(1,0){15}}
\end{picture}}
\newcommand\cellify[1]{\def\thearg{#1}\def\nothing{}%
\ifx\thearg\nothing
\vrule width0pt height\cellsize depth0pt\else
\hbox to 0pt{\usebox2\hss}\fi%
\vbox to 15\unitlength{
\vss
\hbox to 15\unitlength{\hss$#1$\hss}
\vss}}
\newcommand\tableau[1]{\vtop{\let\\=\cr
\setlength\baselineskip{-16000pt}
\setlength\lineskiplimit{16000pt}
\setlength\lineskip{0pt}
\halign{&\cellify{##}\cr#1\crcr}}}
\savebox3{%
\begin{picture}(15,15)
\put(0,0){\line(1,0){15}}
\put(0,0){\line(0,1){15}}
\put(15,0){\line(0,1){15}}
\put(0,15){\line(1,0){15}}
\end{picture}}
\newcommand\expath[1]{%
\hbox to 0pt{\usebox3\hss}%
\vbox to 15\unitlength{
\vss
\hbox to 15\unitlength{\hss$#1$\hss}
\vss}}

\theoremstyle{plain}
\newtheorem{theorem}{Theorem}[section]
\newtheorem{proposition}[theorem]{Proposition}
\newtheorem{prop}[theorem]{Proposition}

\newtheorem{lemma}[theorem]{Lemma}
\newtheorem{cor}[theorem]{Corollary}

\newcommand{\al}{\alpha}
\newcommand{\ual}{\underline{\alpha}}
\newcommand{\oal}{\overline{\alpha}}
\newcommand{\uoal}{\underline{\overline{\alpha}}}
\newcommand{\be}{\beta}
\newcommand{\la}{\lambda}
\newcommand{\ga}{\gamma}

\newcommand{\tga}{\tilde{\gamma}}

\newcommand{\tm}{\tilde{m}}
\newcommand{\lan}{\langle}
\newcommand{\ran}{\rangle}

\newcommand{\epps}{\varepsilon}
\newcommand{\Irr}{\operatorname{Irr}}
\newcommand{\dstyle}{\displaystyle}
\begin{document}

\title[Kronecker products of characters]{On (almost) extreme components in Kronecker products of characters of the symmetric groups}

\author{Christine Bessenrodt}
\address{Institut f\"{u}r Algebra, Zahlentheorie und Diskrete Mathematik,
Leibniz Universit\"{a}t Hannover, Hannover,
D-30167, Germany}
\email{\href{bessen@math.uni-hannover.de}{bessen@math.uni-hannover.de}}

\author{Stephanie van Willigenburg}
\address{Department of Mathematics, University of British Columbia, Vancouver, BC V6T 1Z2, Canada}
\email{\href{mailto:steph@math.ubc.ca}{steph@math.ubc.ca}}

\thanks{The authors thank the Alexander von Humboldt Foundation
and the National Sciences and Engineering Research Council of Canada
for the support of their collaboration.}
\subjclass[2010]{Primary 05E05, 20C30; Secondary 05E10}
\keywords{Kronecker coefficients, Kronecker products,  Schur functions}

\begin{abstract}
{Using a recursion formula due to Dvir, we obtain information
on maximal and almost maximal components in Kronecker products of characters of the
symmetric groups.
This is applied to confirm
a conjecture made by Bessenrodt and Kleshchev in 1999,
which classifies all such Kronecker products
with only three or four components.
}
\end{abstract}

\maketitle

\section{Introduction}\label{sec:intro}

The decomposition of the tensor product of two representations
of a group is an ubiquitous and notoriously difficult problem
which has been investigated for a long time.
For complex representations of a finite group this is equivalent
to decomposing the Kronecker product of their characters
into irreducible characters.
An equivalent way of phrasing this problem for the symmetric groups
is to expand the inner product
of the corresponding Schur functions in the basis of Schur  functions.
Examples for such computations were done
already a long time ago by Murnaghan and Littlewood
(cf.~\cite{Mu,L}).
While the answer in specific cases may be achieved by computing
the scalar product of the characters,
for the important family of the symmetric groups
no reasonable general combinatorial formula is known.

Over many decades, a number of partial results have been
obtained by a number of authors. To name just a few
important cases that come up in the present article,
the products of characters labelled by hook partitions or by two-row partitions
have been computed (see \cite{GR,R,Rosas}), and special constituents, in particular of tensor squares, have been considered \cite{Saxl,V}.
For general products, the largest part and the maximal number of parts
in a partition labelling a
constituent of the product have been determined
(see  \cite{CM,Dvir}); in fact, this is a special case of
Dvir's recursion result \cite[2.3]{Dvir} that
will be crucial in this paper.
The recursion will be used to obtain information on components
to partitions of maximal or almost maximal width and length, respectively.

In general, Kronecker products of irreducible representations
have many irreducible constituents (see e.g. \cite[2.9]{JK}).
In work by Kleshchev and the first author~\cite{BK},
situations are considered where
the Kronecker product of two irreducible $S_n$-characters
has few different constituents.
It was shown there that such products are inhomogeneous
(i.e., they contain at least two different irreducible constituents)
except for the trivial situation where one of the
characters is of degree~1;
indeed, except for this trivial case, no constituent in
a Kronecker product can simultaneously have maximal width and
length.
Investigating the question on homogeneous products
for the representations of the alternating
group $A_n$ motivated the study of products of $S_n$-characters
with two different constituents;
all Kronecker
products of irreducible $S_n$-characters with two homogeneous
components were then classified in~\cite{BK}.
Also, some partial results for products with up to four homogeneous components
were obtained,
and we will use these here.
Moreover, a complete
classification of the pairs $(\chi,\psi)$ of irreducible complex
$S_n$-characters such that the Kronecker product $\chi\cdot\psi$ has three or four homogeneous
components was conjectured in~\cite{BK}.
In this article we
obtain more precise information on the extreme and almost extreme
constituents (and components, respectively) in Kronecker products,
where the {\em extreme constituents} are labelled by partitions of maximal width
or length, while the labels of the {\em almost extreme constituents}  have
width or length differing by one
from the maximal width or length, respectively;
by {\em (almost) extreme constituents (or components)} we mean constituents (or components) of both types.
As a consequence of these results we confirm the conjectures just mentioned.

Note that   $[\mu]$ denotes the irreducible
complex character of $S_n$ labelled by the partition $\mu$ of~$n$;
further details on notation and
background can be found in Section~\ref{sec:prelim}.

{
\begin{theorem}\label{the:extcomp}
Let $\mu$, $\nu$ be partitions of~$n$.
Then the Kronecker product of the characters
$[\mu]$ and $[\nu]$ has at least five
(almost) extreme components unless we are in one
of the following situations.
\begin{enumerate}
\item[(i)] \label{th:1C}
One of $\mu, \nu$ is $(n)$ or $(1^n)$, and the product is irreducible.

\item[(ii)] \label{th:2C}
One of  $\mu, \nu$ is a rectangle $(a^b)$ with $a,b > 1$,
and the other is $(n-1,1)$ or $(2,1^{n-2})$.

\item[(iii)]
$\mu,\nu\in\{(n-1,1),(2,1^{n-2})\}$, for $n>1$.

\item[(iv)]
One of $\mu$, $\nu$ is in $\{(2k,1),(2,1^{2k-1})\}$
while the other one is in $\{(k+1,k),(2^k,1)\}$,
for some $k>1$.

\item[(v)]
Up to conjugating $\mu$ or $\nu$, we have \,
$\{\mu,\nu\} = \{(3^2),(4,2)\}$.

\item[(vi)]
$\mu=\nu$ is a symmetric partition.

\item[(vii)]
$\mu,\nu \in \{((a+1)^a), (a^{a+1})\}$ for some $a>1$.
\end{enumerate}
\end{theorem}

In the cases (i)-(v) above, we know the products explicitly,
while in the last two cases, we know the (almost) extreme components
explicitly.
\smallskip

On the list above, it is not hard to find the products
with exactly three or four components; thus we
obtain the classification conjectured in~\cite{BK}:
}

\begin{theorem}\label{the:34C}
Let $\mu$, $\nu$ be partitions of~$n$.
Then the following holds
for the Kronecker product of the characters
$[\mu]$ and $[\nu]$.
\begin{enumerate}
\item[(i)] \label{th:3C}
The product $[\mu] \cdot [\nu]$ has three homogeneous components
if and only if $n=3$ and $\mu=\nu=(2,1)$ or $n=4$ and $\mu=\nu=(2,2)$.

\item[(ii)]\label{th:4C}
The product $[\mu]\cdot[\nu]$ has four homogeneous components
if and only if $n\geq 4$ and one of the following holds:
\begin{enumerate}[(1)]
\item $\mu,\nu\in\{(n-1,1),(2,1^{n-2})\}$.
\item $n=2k+1$ and
one of $\mu$, $\nu$ is in $\{(n-1,1),(2,1^{n-2})\}$
while the other one is in $\{(k+1,k),(2^k,1)\}$.
\item $\mu,\nu\in\{(2^3),(3^2)\}$.
\end{enumerate}
\end{enumerate}
\end{theorem}

In
Section~\ref{sec:prelim}
we introduce some notation and recall some results
from \cite{BK} and \cite{Dvir} that will be used in the
following sections.
In
Section~\ref{sec:2comp},
we collect results on Kronecker products and
skew characters with at most two components;
on the way, we slightly generalize some of the results
in~\cite{BK} from irreducible characters
to arbitrary characters.
Then
Section~\ref{sec:specprods} deals with special products such as products
with the natural character, squares, products of characters
to 2-part partitions or products of hooks;
this is to some extent based on already available work.
In
Section~\ref{sec:tech},
we put these special cases aside and investigate all other products;
some key results are obtained that help to
produce components in a product which are of
almost maximal width.
These results are then applied
in the proof of Theorem~\ref{the:extcomp},
to which {the main part of} Section~\ref{sec:proof} is devoted.
{Finally, we deduce Theorem~\ref{the:34C}, i.e.,
the conjectured classification
of products with exactly three or four components.}

\section{Preliminaries}\label{sec:prelim}

We denote by $\N$ the set $\{1,2,\dots\}$ of natural numbers,
and let $\N_0=\N \cup\{0\}$.
\smallskip

Let $G$ be a finite group.
We denote by  $\Irr(G)$ the set of irreducible (complex)
characters of~$G$.
Let $\psi$ be a character of $G$.
Then we consider the decomposition of $\psi$
into irreducible characters, i.e.,
$\psi=\sum_{\chi\in \Irr(G)}a_{\chi}\chi$,
with $a_{\chi}\in \N_0$.
If $a_{\chi} >0$, then $\chi$ is called a {\em constituent} of $\psi$
and $a_{\chi}\chi$ is the corresponding  {\em homogeneous component}
of~$\psi$.
The character $\psi$ is called {\em homogeneous} if it has only one homogeneous component.
For any character $\psi$, we let $\X(\psi)$ be the set of irreducible constituents of $\psi$;
so  $c(\psi)=|\X(\psi)|$
is the number of homogeneous components of $\psi$.
If $\psi$ is a virtual character, i.e.,
$\psi=\sum_{\chi\in \Irr(G)}a_{\chi}\chi$
with $a_{\chi}\in \mathbb Z$, then we denote by $c(\psi)$ the
number of components with positive coefficient~$a_{\chi}$.
\smallskip

For the group $G=S_n$, we use the usual notions and notation of
the representation theory of symmetric groups and the related combinatorics
and refer the reader to \cite{JK} or \cite{S} for the
relevant background.
In particular, we write
$\lambda=(\lambda_1, \ldots , \lambda_k) \vdash n$ if
$\lambda$ is a {\em partition} of $n$; in this case we
write $|\lambda|$ for the {\em size} $n$ of~$\lambda$,
$\max \la$ for its largest part or {\em width},
and  $\ell(\lambda)$ for the {\em length} of $\lambda$, i.e.,
the number of its (positive) parts.
We often gather together equal parts of a partition and write
$i^m$ for $m$ occurrences
of the number $i$ as a part of the partition~$\la$.
The partition {\em conjugate} to $\la$ is denoted by $\la'=(\la_1',\ldots,\la'_{\la_1})$.
If $\la=\la'$ we say that $\la$ is {\em symmetric}.
We do not distinguish between a partition $\la$ and its
{\em Young diagram}
$\la=\{(i,j)\in{\N}\times{\N}\mid j\leq\la_i\}$.
Pictorially, we will
draw the diagram using matrix conventions, i.e., starting with row~1 of length $\la_1$ at the top.
Elements $(i,j)\in{\N}\times{\N}$ are called {\em nodes}.
If $\la=(\la_1,\la_2,\dots)$ and $\mu=(\mu_1,\mu_2,\dots)$ are two partitions
we write $\la\cap\mu$ for the partition $(\min(\la_1,\mu_1),\min(\la_2,\mu_2),\dots)$
whose Young diagram is the intersection of the diagrams for $\la$ and $\mu$.
A node
$(i,\la_i)\in\la$ is called a {\em removable} $\la$-node if $\la_i>\la_{i+1}$.
A node
$(i,\la_i+1)$ is called {\em addable} (for $\la$) if $i=1$ or $i>1$ and $\la_i<\la_{i-1}$.
We denote by
$$\la_A=\la\setminus\{A\}=(\la_1,\dots,\la_{i-1},\la_i-1,\la_{i+1},\dots)$$
the partition obtained by removing a removable node  $A=(i,\la_i)$ from $\la$.
Similarly
$$\la^B=\la\cup\{B\}=(\la_1,\dots,\la_{i-1},\la_i+1,\la_{i+1},\dots)$$
is the partition obtained by adding an addable node  $B=(i,\la_i+1)$ to $\la$.

For brevity, when we sum over
all removable $\la$-nodes $A$,
we will sometimes just write
$\dstyle\sum_{A\ \la\text{-node}}$ or $\dstyle\sum_{A}$ to indicate
this summation.

For a node $(i,j)\in \la$,
we denote by
$$
h_{ij}=h_{ij}^\la=\la_i-j+\la_j'-i+1
$$
the corresponding $(i,j)${\em -hook length}.

For a partition $\la$ of $n$,
we write $[\la]$ (or $[\la_1,\la_2,\dots]$) for the (complex) character of $S_n$ associated to $\la$.
Thus, $\{[\la]\mid \la\vdash n\}$ is the set $\Irr(S_n)$ of all (complex)
irreducible characters of $S_n$.

The standard inner product on the class
functions on a group  $G$ is denoted by $\lan\cdot,\cdot\ran$.
If $\chi$ and $\psi$ are two class functions of $G$ we write $\chi\cdot\psi$
for the class function $(g\mapsto\chi(g)\psi(g))$ of $G$.
For characters $[\la]$, $[\mu]$ of $S_n$, the class
function $[\la]\cdot[\mu]$ is again a character of $S_n$, the {\em Kronecker product}
of $[\la]$, $[\mu]$.
We define the numbers
$d(\mu,\nu;\la)$, for $\la,\mu,\nu\vdash n$,   via
$$
[\mu]\cdot[\nu]=\sum_{\la\vdash n} d(\mu,\nu;\la)[\la]\:.
$$

If $\al=(\al_1,\al_2,\dots)$ and $\be=(\be_1,\be_2,\dots)$  are two
partitions then we write $\be\subseteq\al$ if $\be_i\leq \al_i$ for all $i$.
In this case we also consider the skew partition $\al/\be$;
again, we do not distinguish between $\al/\be$ and its skew
Young diagram, which is the set of nodes $\al\setminus\be$
belonging to $\alpha$ but not $\beta$.
If this diagram has the shape of the Young
diagram for a partition, we will also speak of
a {\em partition diagram}.

If $\al/\be$ is a skew Young diagram and $A=(i,j)$ is a node we say $A$ is
{\em connected} to
$\al/\be$ if at least one of the nodes $(i\pm 1,j),(i,j\pm 1)$
belongs to $\al/\be$. Otherwise $A$ is {\em disconnected} from $\al/\be$.

Let  $\be$ and $\ga$ be two partitions. Taking the outer product of the
corresponding characters, the {\em Littlewood-Richardson coefficients}
are defined as the coefficients in the decomposition
$$[\be] \otimes [\ga]=\sum_{\al}c^\al_{\be\ga} [\al]\:.$$
For any partition $\al$,
the {\em skew character} $[\al/\be]$ is then defined to be the sum
$$[\al/\be]=\sum_{\ga}c^\al_{\be\ga} [\ga]\:.$$
The Littlewood-Richardson coefficients can be computed via the
{\em Littlewood--Richardson rule} \cite{JK,M} that says it is the number
of ways of filling the nodes of $\alpha/\beta$ with positive integers
such that the rows weakly increase left to right, the columns strictly
increase top to bottom, and when the entries are read from right to left
along the rows starting at the top, the numbers of $i$'s read is always
weakly greater than the number of $(i+1)$'s.
Note that $[\al/\be]=0$ unless $\be\subseteq\al$.
\\

As in \cite{BK}, we will repeatedly use the following results.
The first describes the rectangular hull of the partition
labels of the constituents of a Kronecker product, i.e.,
the \emph{maximal width} and \emph{maximal length} of these labels.

\begin{theorem}\label{ThDvir1} \cite[1.6]{Dvir}, \cite[1.1]{CM}.
Let $\mu$, $\nu$ be partitions of $n$.
Then
$$\max\{\max\la \mid  [\la ]\in \X([\mu]\cdot[\nu]) \}
    = |\mu \cap \nu|$$
and
$$\max\{\ell(\la) \mid   [\la ]\in \X([\mu]\cdot[\nu])  \}
    = |\mu \cap \nu'|\:.$$
\end{theorem}

If a component or constituent is described as being of \emph{almost} maximal width or length the respective width or length differs from the above values by one.

Since skew characters of $S_n$ can  be decomposed
into  irreducible characters using the Littlewood-Richardson rule,
the following theorem provides a
recursive formula for the coefficients $d(\mu,\nu;\la)$.

{
In the following, if $\la=(\la_1,\la_2,\dots,\la_m)$ is a partition of~$n$,
of length~$m$,
we set $\hat\la = (\la_2, \la_3 , \ldots)$
and $\bar{\la}=(\la_1-1,\la_2-1,\dots,\la_m-1)$.
}

\begin{theorem}\label{ThDvir2} \cite[2.3]{Dvir}.
Let $\mu$, $\nu$ and $\la=(\la_1,\la_2,\dots)$ be partitions of~$n$.
Define
$$Y(\la) = \{ \eta=(\eta_1,\ldots) \vdash n \mid
  \eta_i \geq \la_{i+1} \geq \eta_{i+1} \mbox{ for all } i \geq 1\}\:,$$
i.e., $Y(\la)$ is the set of partitions obtained from $\hat\la$ by adding a horizontal strip of size~$\la_1$.
Then
$$ d(\mu,\nu;\la)
= \sum_{{\alpha \vdash \la_1}\atop {\alpha \subseteq \mu \cap \nu}}
 \langle [\mu / \alpha] \cdot [\nu/\alpha] , [\hat{\la}]\rangle
  - \sum_{{{\eta \in Y(\la) \atop {\eta \neq \la}} \atop
  \eta_1 \leq |\mu \cap \nu|} }  d(\mu,\nu;\eta) \; .$$
\end{theorem}

{
\begin{cor}
Let $\mu,\ \nu,\ \la$ be partitions of~$n$, and set
$\gamma = \mu \cap \nu$.
\begin{enumerate}
\item[(i)] \label{CDvir2} (\cite[2.4]{Dvir}, \cite[2.1(d)]{CM})
If  $\max\la = |\mu \cap \nu|$, then
$$ d(\mu,\nu;\la)
= \langle [\mu / \gamma] \cdot [\nu / \gamma] , [\hat{\la}]\rangle
\; .$$
\item[(ii)]\label{CDvir2'} (\cite[2.4']{Dvir})
If $\ell(\la)=|\mu \cap \nu'|$, then
$$
d(\mu,\nu;\la)=\langle [\mu/(\mu\cap\nu')]\cdot[\nu/(\mu'\cap\nu)],[\bar{\la}]\rangle\:.
$$
\end{enumerate}
\end{cor}
}
\smallskip

We will later use these results for finding {\em extreme constituents}
in a product, i.e., those with partition labels of  maximal width or  maximal length;
explicitly, we state this in the following lemma.
Here, for a partition $\al=(\al_1,\al_2,\ldots)$ and $m\in \N$,
we set {$(m,\al)=(m,\al_1,\al_2,\ldots)$ and}
$\al + (1^m) = (\al_1+1, \ldots , \al_m+1, \al_{m+1}, \ldots)$.
\begin{lemma}\label{extcomp}
Let $\mu,\ \nu,\ \la$
be partitions of~$n$, $\ga=\mu \cap \nu \vdash m$,
$\tga=\mu \cap \nu' \vdash \tm$.
\begin{enumerate}[(i)]
\item If $[\alpha]$ appears in $[\mu / \ga] \cdot [\nu / \ga]$,
then $[m,\al]$  appears in $[\mu]\cdot[\nu]$.
\item If $[\alpha]$ appears in $[\mu / \tga]  \cdot [\nu' / \tga]$,
then $[\al' + (1^{\tm})]$  appears in $[\mu]\cdot[\nu]$.
\end{enumerate}
\end{lemma}

\proof
(i) If  $[\alpha]$ appears in $[\mu/\ga]\cdot [\nu/\ga]$, then
it is a constituent of some $[\rho]\cdot [\tau]$ with
$\rho \subset \mu$, $\tau\subset \nu$
by the Littlewood-Richardson rule,
and hence by Theorem~\ref{ThDvir1}
$\al_1 \leq |\rho\cap\tau| \leq |\mu\cap\nu|=m$.
Thus $(m,\al)$ is a partition and
the claim follows by Corollary~\ref{CDvir2}(i).

(ii)  If $[\alpha]$ appears in $[\mu / \tga]  \cdot [\nu' / \tga]$,
then by Theorem~\ref{ThDvir1}
$\al_1 \leq  |\mu\cap\nu'|=\tm$.
Hence $\al' + (1^{\tm})$ is a partition and
the claim follows by Corollary~\ref{CDvir2'}(ii).
\qed

\medskip

The following result implies the important fact
that in {\em nontrivial}
Kronecker products $[\mu]\cdot[\nu]$,
i.e., products where both characters are of degree $>1$,
there is no constituent which is both of maximal width and maximal length:
\begin{theorem}\label{PConst}\cite[Theorem 3.3]{BK}
Let $\mu$, $\nu$ be partitions of~$n$, both different from
$(n)$ and~$(1^n)$.
If $[\lambda]$ is a constituent of $[\mu] \cdot [\nu]$,
then $h_{11}^{\lambda} <     |\mu \cap \nu| + |\mu \cap \nu'| - 1$.
\end{theorem}

The constituents of maximal width (or length) in a Kronecker product may be handled
by the results above.

We also note the following consequence of
Theorem~\ref{ThDvir2} which indicates how
to find components of almost maximal width.

\begin{lemma}\label{lem:prepalmext}
Let $\mu, \nu , \la=(\la_1, \ldots ) \vdash n$,
$\ga=\mu \cap \nu \vdash m$; assume $\la_1=m-1$.
We define the virtual character
$$
\chi=
\sum_{A\ \mathrm{ removable}\ \ga\mathrm{-node }} [\mu / \ga_A] \cdot [\nu / \ga_A]
- ([\mu/\ga]\cdot [\nu/\ga])\uparrow^{S_{n-m+1}}
\:.
$$
Then
\; $ d(\mu,\nu;\la) =  \langle\chi,[\hat\la]\rangle$.
\end{lemma}

\proof
For the recursion formula in Theorem~\ref{ThDvir2} we need to compute
$\dstyle \epps=   \sum_{\eta \in \tilde Y (\la) }  d(\mu,\nu;\eta)$,
where \;
$\dstyle
\tilde Y(\la) = \{  \eta \vdash n \mid
  \eta_i \geq \la_{i+1} \geq \eta_{i+1} \mbox{ for all } i \geq 1,
  \eta_1\leq m, \eta\neq \lambda \}$.

By definition,
$\eta\in \tilde Y(\la)$ arises from $\hat\la$ by putting on a row strip
of size $\la_1=m-1$, but as $\eta_1\leq m$,
we must then have $\eta=(m,\hat\la_A)$ for some removable $\hat\la$-node~$A$.

Let \;
$\dstyle [\mu/\ga] \cdot [\nu/\ga]= \sum_{\al} c_{\al} [\al]$;
then by Corollary~\ref{CDvir2}(i)
$$[\mu]\cdot[\nu]= \sum_{\al} c_{\al} [m,\al] +
\sum_{\be \vdash n \atop {\beta_1<m}} d(\mu,\nu;\be)\ [\be] \:.
$$
Thus
$$\begin{array}{rcl}
\epps &=&
\dstyle
\sum_{\eta \in \tilde Y (\la) }  d(\mu,\nu;\eta)
=
\sum_{\eta \in \tilde Y (\la) }  \sum_{\al} c_{\al} \langle[m,\al], [\eta]\rangle
=
\sum_{A \, \hat\la\text{-node} }  \sum_{\al} c_{\al} \langle[\al], [\hat\la_A]\rangle
\\[20pt]
&=&
\dstyle
\langle [\mu/\ga]\cdot[\nu/\ga], \sum_{A \, \hat\la\text{-node} } [\hat\la_A]\rangle
= \langle ([\mu/\ga]\cdot[\nu/\ga])\uparrow^{S_{n-m+1}}, [\hat\la]\rangle \:.
\end{array}
$$
As $\la_1=m-1$, the assertion on $d(\mu,\nu;\la)$ now follows
from Theorem~\ref{ThDvir2}.
\qed
\medskip

In \cite{BK}, some information on constituents  of almost maximal width
was obtained in the special situation described below
which will be used here as well.
Beware that here we do not start with partitions $\la$ corresponding
to (almost) maximal width constituents in $[\mu]\cdot [\nu]$,
but the subtle point  in the statement below is that $(m,\al)$ and $(m-1,\theta)$
are indeed partitions.

\begin{lemma}\label{BK4.6}\cite[Lemma 4.6]{BK}
Let $\mu \neq \nu$ be partitions of~$n$,
both different from   $(n)$, $(1^n)$, $(n-1,1)$ and $(2,1^{n-2})$.
Put $\gamma = \mu \cap \nu$, $m = |\gamma|$.
Assume that $\nu / \gamma$ is a row and that $[\mu / \gamma]$
is an irreducible character $[\al]$.
Then $[m,\al] \in \X([\mu]\cdot[\nu])$.

Moreover, if an $S_{n-m+1}$-character $[\theta]$ appears in
\begin{equation}\label{virtual}
\chi = 
\sum_{A\ \mathrm{ removable}\ \ga\mathrm{-node} }[\mu / \gamma_A] \cdot [\nu / \gamma_A]
- [\al]\uparrow^{S_{m-n+1}}
\end{equation}
with a positive coefficient, then $[m-1,\theta]\in \X([\mu]\cdot[\nu])$.
\end{lemma}

\section{Products and skew characters with at most two components}
\label{sec:2comp}

In \cite{BK}, the Kronecker products of $S_n$-characters with at most two components
were classified;
in the nontrivial case, then exactly one component is of maximal width
and one {is} of maximal length.
We recall the classification here explicitly
as we will need this {later}.

\begin{theorem}\label{Th1C}\cite[Corollary 3.5]{BK}
Let $\mu$, $\nu$ be partitions of~$n$.
Then $[\mu] \cdot [\nu]$ is homogeneous if and only if
one of the partitions is $(n)$ or $(1^n)$.
In this case $[\mu] \cdot [\nu]$ is irreducible.
\end{theorem}

\begin{theorem}\label{Th2C}\cite[Theorem 4.8]{BK}
Let $\mu$, $\nu$ be partitions of~$n$.
Then $[\mu] \cdot [\nu]$ has exactly two homogenous components
if and only if
one of the partitions $\mu, \nu$ is a rectangle $(a^b)$ with $a,b > 1$,
and the other is $(n-1,1)$ or $(2,1^{n-2})$.
In these cases we have:
\begin{eqnarray*}
& & [n-1,1] \cdot [a^b] = [a+1,a^{b-2},a-1] + [a^{b-1},a-1,1] \\
& & [2,1^{n-2}] \cdot [a^b] = [b+1,b^{a-2},b-1] + [b^{a-1},b-1,1]\:.
\end{eqnarray*}
\end{theorem}

\medskip

Towards our understanding of (almost)
extreme components, we need more information
on skew characters.
We consider here in particular the situation of
skew characters with few components.
>From \cite{BK} we quote the following (see also \cite{vW_schur}):

\begin{lemma}\label{LSkCh}\cite[Lemma 4.4]{BK}
Let $\mu$, $\ga$ be partitions, $\ga \subset \mu$.
Then the following assertions are equivalent:
\begin{enumerate}[(a)]
\item $[\mu / \ga]$ is homogeneous.
\item $[\mu / \ga]$ is irreducible.
\item The skew diagram $\mu / \ga$ is the diagram of
a partition $\al$ or the rotation of a diagram of a partition $\al$.
In this case,  $[\mu / \ga] = [\al]$.
\end{enumerate}
\end{lemma}

\medskip

Note that throughout this paper, whenever we say rotation we mean
{\em rotation by $180^\circ$}.
For the later consideration of almost extreme components,
we will also need the following from~\cite{BK}.

\begin{lemma}\label{lem:4.5'}
Let $\mu$, $\ga$ be partitions, $\ga \subset \mu$,
such that $[\mu/\ga]$ is irreducible, say $[\mu/\ga]=[\al]$.
Let $A$ be a removable node of $\ga$.
\begin{enumerate}[ (1)]
\item
If $A$ is disconnected from $\mu / \ga$ then
$$
[\mu / \ga_A] =  \sum_{B\ \mathrm{addable\ for }\ \al} [\al^B]\:.
$$
\item Let $A$ be connected to $\mu / \ga$.
Let $B_0$ and $B_1$ denote the top and bottom removable node of $\al$
 respectively.

If $\mu/\ga$ has partition shape  $\al$, then
$$[\mu / \ga_A] =
\left\{
\begin{array}{ll}
\dstyle \sum_{B\ \mathrm{addable\ for}\ \al \atop {B \neq B_0}} [\al^B] &
\text{if $A$ is connected to the top row of $\al$}\\[20pt]
\dstyle \sum_{B\ \mathrm{addable\ for}\ \al \atop {B \neq B_1}} [\al^B] &
\text{if $A$ is connected to the bottom row of $\al$.}\\
\end{array}
\right.
$$

If $\mu/\ga$ is the rotation of $\al$, then
$$[\mu / \ga_A] =  [\al^B] \:, \:
\text{where $B$ is an addable node of $\al$.}$$
\end{enumerate}
\end{lemma}

\medskip

{As seen above, the investigation of extreme components
in Kronecker products  requires
a study of the components of skew characters, and we now turn to the case
where these have two components.}

We denote an outer product of two characters $\chi$, $\psi$ by $\chi \otimes \psi$; recall that an outer product of two irreducible characters corresponds to a character associated to a skew diagram decomposing into two
disconnected partition diagrams.
In the classification list below,
$r,s,a,b$ are arbitrary nonnegative integers
such that all characters appearing on the left hand side correspond to partitions.
For example,  $[(r+1)^{a+1},r^{b},s]+[(r+1)^{a},r^{b+1},s+1]=[((r+1)^{a+1},r^{b+1})/(r-s)]$
also incorporates that $r\geq s+1$ since otherwise $((r+1)^{a},r^{b+1},s+1)$ is not a partition; thus here $s,a,b\geq 0$ and $r\geq s+1$.
Choosing the minimal values $r=1,s=a=b=0$ gives $[2]+[1^2]=[(2,1)/(1)]$.
Note that any proper skew character $[\mu/\ga]$, i.e., one that is not irreducible,
always has a constituent obtained by sorting the row lengths of $\mu/\ga$
as well as a (different) constituent obtained by sorting the column lengths and
conjugating{, and both are of multiplicity~1}.
These appear on the right hand side below, which we
write out explicitly for the convenience of the reader.
{From the work of Gutschwager we have the following
classification list for skew characters with exactly two components.}

\begin{prop}\label{SkCh2Comp}\cite{CG}
The following is a complete list of skew characters of symmetric groups
with  exactly two homogeneous components; up to rotation, ordering and translation
all corresponding skew diagrams are given.
 \begin{enumerate}[(i)]
  \item $[1]\otimes[r^{a+1}]=[r+1,r^a]+[r^{a+1},1]$.
  \item  $[1^{a+1}]\otimes[r]=[r,1^{a+1}]+[r+1,1^a]$.
  \item  $[((r+1)^{a+1},r^{b+1})/(r-s)]= [((r+1)^{a+b+1},s+1)/(1^{b+1})]$ \\ $=[(r+1)^{a+1},r^{b},s]+[(r+1)^{a},r^{b+1},s+1]$.
  \item  $[((r+1)^{a+1},(s+1)^{b+1})/(1)]=
      [(r+1)^a,r,(s+1)^{b+1}]+[(r+1)^{a+1},(s+1)^b,s]$.
  \item  $[((r+1)^{a+1},(s+1)^{b+1})/(1^{a+b+1})]=[r^{a+1},s+1,s^b]+[r+1,r^a,s^{b+1}]$.
  \item $[(r^{a+1},s^{b+1})/(r-1)]=[r^a,s^{b+1},1]+[r^a,s+1,s^b]$.
  \item $[(r+1,(s+1)^{a+b+1})/(1^{b+1})]=[r,(s+1)^{a+1},s^b]+[r+1,(s+1)^a,s^{b+1}]$.
  \item  $[(r^{a+1},1^{b+1})/(r-s)]=[r^a,s,1^{b+1}]+[r^a,s+1,1^b]$.
 \end{enumerate}
\end{prop}
\medskip

For later purposes we note the following slight generalization
of Theorem~\ref{Th1C};
recall that $c(\chi)$ denotes the number of homogeneous components
of a character $\chi$.

\begin{prop}\label{Prop:homocharprod}
Let $\chi$, $\psi$ be characters of~$S_n$.
Then $c(\chi \cdot \psi)=1$  if and only if
we are in one of the following situations (up to
ordering the characters):
\begin{enumerate}[(1)]
\item $c(\chi)=1$ and $\psi=a[n]$ or $\psi=b[1^n]$, with $a,b \in \N$.
\item $\chi=k[\alpha]$ with $\alpha=\alpha'$, $k\in \N$,
and $\psi=a[n]+b[1^n]$, with $a,b\in \N_0$.
\end{enumerate}
\end{prop}

\proof
This follows easily using Theorem~\ref{Th1C}.
\qed

\begin{cor}\label{Cor:skewdia1C}
Let $D, \tilde{D}$ be skew diagrams of size~$n$.
Then $c([D]\cdot [\tilde{D}])=1$ if and only if
one of the diagrams is a partition or rotated partition diagram,
and the other one is a row or column.

In particular, the product of two skew characters is
a homogeneous character
if and only if it is an irreducible character.
\end{cor}

\proof
Use Proposition~\ref{Prop:homocharprod} and Lemma~\ref{LSkCh}, and
note that a skew character can be of the form $a[n]+b[1^n]$ with $a,b>0$
only for $n=2$, where there is no symmetric partition.
\qed

\medskip

Using Theorem~\ref{Th2C}, one may generalize the result and
characterize the products of characters with
exactly two homogeneous components:

\begin{theorem}\label{Thm:2cp}
Let $\chi,\psi$ be characters of $S_n$, $n>1$.
Then $c(\chi \cdot \psi)=2$ exactly in the following situations
(up to ordering the characters):
\begin{enumerate}
\item
$\chi=a[\la]$ with $\la\neq \la'$, $\psi=b[n]+c[1^n]$,
$a,b,c\in \N$.
\item
$\chi=a[r^s]$ with $r\neq s$, $r,s>1$, $\psi=b[n-1,1]$ or $b[2,1^{n-2}]$,
$a,b\in \N$.
\item
$\chi=a[r^r]$, $\psi=b[n-1,1]+c[2,1^{n-2}]$, $a,b,c \in \N_0$, $a, b+c>0$.
\item
$\chi=a[n]$ or $a[1^n]$, $\psi=b[\la]+c[\mu]$, $\la \neq \mu$ arbitrary,
$a,b,c\in \N$.
\item
$\chi=a[n]+b[1^n]$, $\psi=c[\la]+d[\la']$ with $\la\neq \la'$,
$a,b,c,d\in \N$.
\item
$\chi=a[n]+b[1^n]$, $\psi=c[\alpha]+d[\beta]$ with $\alpha=\alpha'\neq \beta=\beta'$, $a,b,c,d \in \N$.
\item
$\chi=a[2^2]+b[4]+c[1^4]$,  $\psi=d[3,1]+e[2,1^2]$,
$a,b,c,d,e\in \N_0$, $a, b+c,d+e >0$.
\end{enumerate}
\end{theorem}

\medskip

In the following, by a {\em nontrivial rectangle}  we mean a rectangle
with at least two rows and at least two columns.

\begin{cor}\label{Cor:skewdia2C}
Let $D, \tilde{D}$ be skew diagrams of size~$n$.
Then $c([D]\cdot [\tilde{D}])=2$ if and only if
we are in one of the following situations (up to ordering
of the diagrams):
\begin{enumerate}
\item
$n=2$, $[D]=[\tilde{D}]=[2]+[1^2]$; here, $D$ and $\tilde{D}$
both consist of two disconnected nodes.
\item
$D$ is a row or column, $[\tilde{D}]$ is one of the skew characters in
Proposition~\ref{SkCh2Comp}; the diagram
$\tilde{D}$ corresponds to one of the diagrams in
Proposition~\ref{SkCh2Comp} up to rotation,
translation and reordering disconnected parts.
\item
$D$ is a nontrivial rectangle, $\tilde{D}$ has shape $(n-1,1)$ up to rotation and conjugation.
\item
$n=4$, $D=(2^2)$, $\tilde{D}$ consists of a disconnected row and column of size~2 each; here, $[\tilde{D}]=[3,1]+[2,1^2]$.
\end{enumerate}
\end{cor}

\proof
Theorem~\ref{Thm:2cp} strongly restricts the number of cases to be
considered.
Then, note that a skew character $[D]$
can be of the form $a[n]+b[1^n]$ with $a,b>0$ only
for $n=2$, when  $D$ consists of two disconnected nodes.
In fact, if both $[n]$ and $[1^n]$ appear in the skew character,
then all irreducible characters appear as constituents.
The only case where both  $[D]$ and $[\tilde D]$ are inhomogeneous is
the one described for $n=2$ in~(1) above.
Homogeneous skew characters
are labelled by partition diagrams, up to rotation,
and are indeed irreducible (by Lemma~\ref{LSkCh}).
Furthermore, we use Proposition~\ref{SkCh2Comp} to identify
the diagrams for the skew characters with two components
appearing above, and we note that $a[2^2]+b[4]$ and $a[2^2]+c[1^4]$ are
not  skew characters.
\qed

\section{Special products}\label{sec:specprods}

The products  with the character $[n-1,1]$
are easy to compute, and
{those products with few components have been
classified \cite{BK}; we extend this here a bit further.}

\begin{lemma}\cite[Lemma 4.1]{BK}\label{LNat}
Let $n \geq 3$ and $\mu$ be a partition of~$n$.
Then
$$[\mu] \cdot [n-1,1] = \sum_A \sum_B [(\mu_A)^B] - [\mu]$$
where the first sum is over all removable nodes $A$ for $\mu$,
and the second sum runs over all addable nodes~$B$ for $\mu_A$.
\end{lemma}

\begin{cor} (see \cite[Cor.~4.2]{BK})\label{Cor:Nat}
Let $n \geq 3$ and $\mu$ be a partition of~$n$.
{If $[n-1,1] \cdot [\mu]$ has at most four (almost) extreme
components then the total number of components is at most four
and we are in one of the following situations.
}
\begin{enumerate}[(i)]
\item $c([\mu] \cdot [n-1,1])=1$
if and only if $\mu$ is $(n)$ or $(1^n)$.
\item $c([\mu] \cdot [n-1,1])=2$
if and only if $\mu$ is a rectangle
$(a^b)$ for some $a,b > 1$. In this case we have
$$[a^b] \cdot [n-1,1] =
[a+1,a^{b-2},a-1] + [a^{b-1}, a-1,1]\:.$$
\item $c([\mu] \cdot [n-1,1])=3$
if and only if $n=3$ and $\mu = (2,1)$.
In this case we have
$$[2,1] \cdot [2,1] = [3] + [2,1] + [1^3]\:.$$
\item $c([\mu] \cdot [n-1,1])=4$
if and only if one of the following happens:
\begin{enumerate}[(a)]
\item $n \geq 4$ and $\mu = (n-1,1)$ or $(2, 1^{n-2})$.
\item $\mu = (k+1,k)$ or $(2^k,1)$ for $k \geq 2$.
\end{enumerate}

We then have:
$$\begin{array}{rcl}
[n-1,1] \cdot [n-1,1] & = &
[n]+[n-1,1]+[n-2,2]+[n-2,1^2]\\[5pt]
[k+1,k] \cdot [2k,1]  & =  & [k+2,k-1]+[k+1,k]+[k+1,k-1,1]+[k^2,1]
\end{array}$$
and the remaining products are obtained by conjugation.
\end{enumerate}
{
(Note that in all cases above all components are (almost) extreme.)
}
\end{cor}

{
\proof
We may assume that we are not in one of the situations (i)-(iv) described above,
and we want to show that we then have at least five (almost) extreme components.
The extreme components of $[\mu] \cdot [n-1,1]$
have width $\mu_1+1$ and length $\ell(\mu)+1$, respectively.
Let $X$ and $Y$ be the top and bottom addable nodes for~$\mu$, and $r$ the number
of removable nodes of~$\mu$; note that $r>1$ as $\mu$ is not a rectangle, by our assumption.
With $A$ running over the removable nodes of~$\mu$, but not in the first row, we get $r'\in \{r,r-1\}$ constituents $[(\mu_A)^X]$ of maximal width in the product;
similarly, with $B$ running over the removable nodes of~$\mu$, but not in the first column,
we get $r''\in \{r,r-1\}$ constituents $[(\mu_B)^Y]$ of maximal length in the product.
Note that $[\mu]$ itself appears as an almost extreme component in the product,
with multiplicity $r-1\geq 1$.
Thus, the only critical cases to be discussed are the ones where $r'+r''\in \{2,3\}$.

The case $r'+r''=2$ occurs only for hook partitions $\mu$, say $\mu=(n-k,1^k)$,
with $1<k<n-2$, by our assumption.
Here, the product is
$$
[n-k,1^k]\cdot [n-1,1] = [n-k+1,1^{k-1}]+[n-k,2,1^{k-2}]+[n-k,1^k]+[n-k-1,2,1^{k-1}]
+[n-k-1,1^{k+1}]\: ,
$$
with five components which are all (almost) extreme.

The case $r'+r''=3$ occurs only for $r=2$, i.e., proper fat hooks $\mu$,
where either $\mu_1=\mu_2$ but {${\mu '}_1>{\mu '}_2$ or  $\mu_1>\mu_2$ and ${\mu '}_1={\mu '}_2$.}
Using conjugation, we may assume that we are in the first case; then
$\mu=(a^b,1^c)$ for parameters $a,b>1$, $c\geq 1$, by our assumption.
Let $Z$ be the middle addable node, and $A$, $B$ the top and bottom removable
nodes for $\mu$, respectively. If $c>1$, then $[(\mu_B)^Z]$ is a component
of the product of almost maximal width; if $c=1$, then $a>2$, by our assumption
that we are not in situation~(iv),
and  $[(\mu_A)^Z]$ is a component of almost maximal width.
Hence in any case we have found a fifth (almost) extreme
component.

Thus, in any situation different from the ones described in (i)-(iv),
we have found five (almost) extreme
components in the product.
\qed
}

\medskip

As a further special family of products the Kronecker squares were
considered in \cite{BK}{, and the ones with at most four components
were already classified}:
\begin{proposition}\cite[Lemma 4.3]{BK} \label{LTensSq}
Let $\lambda$ be a partition of~$n$.
Then $c([\lambda]^2)\leq 4$
if and only if one of the following holds:
\begin{enumerate}[(1)]
\item $\lambda = (n)$ or $(1^n)$, when $[\lambda]^2=[n]$.
\item $n \geq 4$, $\lambda = (n-1,1)$ or $(2,1^{n-2})$,
            when $[\lambda]^2 = [n]+[n-1,1]+[n-2,2]+[n-2,1^2]$.
\item $n=3$, $\lambda = (2,1)$,
          when $[\lambda]^2 = [3]+[2,1]+[1^3]$.
\item $n=4$, $\lambda = (2^2)$,
           when $[\lambda]^2 = [4]+[2^2]+[1^4]$.
\item  $n=6$, $\lambda = (3^2)$ or $(2^3)$,
           when $[\lambda]^2 = [6]+[4,2]+[3,1^3]+[2^3]$.
\end{enumerate}
\end{proposition}

{Note that here, the squares of $[2^2]$ and $[3^2]$ do have
components which are not (almost) extreme.}

We need to extend the result above further for its
application in the proof of Theorem~\ref{the:extcomp}:
\begin{proposition}\label{prop:TensSqEx}
Let $\lambda$ be a partition of~$n$.
Then $[\lambda]^2$ has at most four (almost) extreme components
if and only if one of the following holds:
\begin{enumerate}[(1)]
\item $\lambda = (n)$ or $(1^n)$.
\item $n \geq 4$, $\lambda = (n-1,1)$ or $(2,1^{n-2})$.
\item $n>1$, $\lambda = \lambda'$.
\item $a>1$, $\lambda = (a^{a+1})$ or $((a+1)^a)$.
\end{enumerate}
\end{proposition}

\proof
We may assume that $\la$ is not one of the partitions already appearing in the
list of Proposition~\ref{LTensSq}.

First, let $\la=\la'$, i.e., $\la$ is a symmetric partition, with $r$ removable nodes.
Using Lemma~\ref{LNat}  and symmetry of Kronecker coefficients we have
$$[\la]^2 = [n] + (r-1)[n-1,1]+ \text{(other components)} +
(r-1)[21^{n-2}]+[1^n]\:,$$
and clearly, the other components appearing here are not (almost) extreme,
so we have at most four (almost) extreme components.

\smallskip

Now let $\la$ be a nonsymmetric partition with $r$ removable nodes. Then
$$[\la]^2 = [n] + (r-1)[n-1,1]+ \text{(other components)} \:.$$
The components of (almost) maximal length in $[\la]^2$ come from considering
the components of (almost) maximal width  in
$[\la]\cdot [\la']$.
Using Lemma~\ref{extcomp}
we get the ones of maximal width $m=|\la\cap \la'|$
from the components of
$[\la\setminus (\la \cap \la')]\cdot [\la'\setminus (\la \cap \la')]$.
This  product has almost always at least three components;  by
Corollaries~\ref{Cor:skewdia1C} and~\ref{Cor:skewdia2C}
and Theorem~\ref{Thm:2cp} the only
cases with at most two components are of the following types (up to conjugation):

(i) $\la\setminus (\la \cap \la')$ is a row;

(ii) $\la\setminus (\la \cap \la')$ consists of two disjoint nodes.

Note that by our assumption $m\geq 4$, so the (almost) maximal width components
in $[\la]\cdot [\la']$ give (almost) maximal length components in $[\la]^2$
different from the (almost) maximal width components already given above.
We get only one or two maximal width components in the product  $[\la]\cdot [\la']$
in cases (i) and (ii) above, and we need to find further almost maximal
width components.

\smallskip

{\bf Case (i).} $\la\setminus (\la \cap \la')$ is a row.

Set $\ga=\la\cap \la' \vdash m$.
With $[\al]=[\la'/\ga]=[1^{n-m}]$, we can apply Lemma~\ref{BK4.6}:
components of
$$\chi=\sum_A [\la/\ga_A]\cdot [\la'/\ga_A] - [21^{n-m+1}]-[1^{n-m+1}]$$
induce components of almost maximal width in $[\la]\cdot [\la']$.

Assume that there is a removable $\ga$-node $A$
not connected to $\la/\ga$ and $\la'/\ga$.
Then we have
$$
\begin{array}{rcl}
[\la/\ga_A]\cdot [\la'/\ga_A]
& = & ([n-m+1]+[n-m,1])\cdot([1^{n-m+1}]+[21^{n-m-1}])
\\[5pt]
& = & 2[1^{n-m+1}]+3[21^{n-m-1}]+[31^{n-m-2}]+[2^21^{n-m-3}]
\end{array}
$$
where the third and fourth component only appear for $n-m\geq 2$ and $n-m\geq 3$,
respectively.
As $n-m\geq 1$, we thus find at least two components in $\chi$;
hence we have at least three (almost) maximal components in $[\la]\cdot [\la']$.
As $\la$ cannot be a rectangle {(otherwise $\ga$ would be a square, with the only
corner node connected to $\la/\ga$)}, this gives altogether at least five (almost) extreme
components in the square $[\la]^2$.
\smallskip

We may now assume that all removable $\ga$-nodes are connected to $\la/\ga$ or
$\la'/\ga$ (or both).

Assume that there is a removable $\ga$-node $A_0$
connected to $\la/\ga$, but not to
$\la'/\ga$; by the symmetry of $\ga$, we then
have a conjugate $\ga$-node $A_1$ connected to $\la'/\ga$, but not to
$\la/\ga$.

If $\la/\ga_{A_0}$ still is a row, we have
$$
[\la/\ga_{A_0}]\cdot [\la'/\ga_{A_0}] + [\la/\ga_{A_1}]\cdot [\la'/\ga_{A_1}]
= 2([1^{n-m+1}]+[21^{n-m-1}])\:,
$$
and thus $\chi$ has at least two components, and we are done as before.
\smallskip

We are left with the case where all removable $\ga$-nodes are connected to
both $\la/\ga$ and~$\la'/\ga$. Clearly, this can only happen when there is only one
such $\ga$-node, $A$ say, that is the corner of the square~$\ga =(a^a)$,
and then $\la=(a^{a+1})$; by our assumption, $a>2$ and hence $n-m>2$.
In this situation we have
$$
\chi = [n-m,1]\cdot [2 1^{n-m-1}] - [2 1^{n-m-1}] - [1^{n-m+1}]
= [3,1^{n-m-2}] + [2^2 1^{n-m-3}]
$$
and hence we have indeed exactly three components of (almost) maximal width
in $[\la]\cdot [\la']$ and thus three components of (almost) maximal length
in  the square $[\la]^2$.
As $\la$ is a rectangle, we have only one component of (almost) maximal width in the square,
hence altogether exactly four (almost) extreme components.

\smallskip

{\bf Case (ii).}  $\la\setminus (\la \cap \la')$ consists of two disjoint nodes.
{Then by Corollary~\ref{CDvir2}(i)}
$$[\la] \cdot [\la'] = 2[m,2]+2[m,1^2] + \text{(further components)} \:,
$$
so we have two components in $[\la]^2$ of maximal length~$m \geq 4$.
As $\ga$ has at least four addable nodes, it has at least three removable nodes.

Assume first that there is a removable $\ga$-node $A_0$ not connected to $\la/\ga$.
Then $[\la/\ga_{A_0}]\cdot [\la'/\ga_{A_0}]$ contains the subcharacter
$$
([3]+2[2,1]+[1^3])\cdot [2,1] = 2[3]+4[2,1]+2[1^3] = ([\la/\ga]\cdot[\la'/\ga])\uparrow^{S_3}\:.
$$
For any other removable $\ga$-node $A_1$, $[\la/\ga_{A_1}]\cdot [\la'/\ga_{A_1}]$
contains the subcharacter $[2,1]^2=[3]+[2,1]+[1^3]$.
Hence the character $\chi$ as defined in Lemma~\ref{lem:prepalmext}
has three components, inducing three components
of almost maximal width $m-1\geq 3$ in $[\la] \cdot [\la']$
and hence three components
of almost maximal length $m-1 \geq 3$ in $[\la]^2$, and we are done.

If there is no such $\ga$-node $A_0$, for both $\la/\ga$ and $\la'/\ga$,
then $\ga=(3,2,1)$ and $\{\la,\la'\}=\{(4,2^2), (3^2,1^2)\}$.
In this case we have
$$
\chi = 3([3]+[2,1])\cdot ([1^3]+[2,1])-2([2]+[1^2])\uparrow^{S_3}
= [3]+5[2,1]+4[1^3]\:,
$$
and using Lemma~\ref{lem:prepalmext} again and conjugation,
we get three components of almost maximal length~5 in $[\la]^2$,
finishing the proof.
\qed
\medskip

We now turn to further families of Kronecker products where we can classify the
products with few components.
While the products of characters to 2-part partitions and hooks have been determined
in work by Remmel et al. \cite{R, R92, RW}
and Rosas~\cite{Rosas},
here we do not use these intricate results but prove the following weaker
facts for the sake of a self-contained presentation.
\smallskip

By a {\em 2-part partition} we mean here a partition of length at most~2;
we say it is {\em proper} if it has length exactly~2.

\begin{proposition}\label{twopart}
Let $\mu, \nu \vdash n$ be different proper 2-part partitions,
say $\mu=(n-k,k)$, $\nu=(n-l,l)$ with $1<l<k$.
Let $\ga=\mu\cap\nu$, a partition of $m=n-k+l$.
Then we have the following constituents in $[\mu]\cdot[\nu]$:
\begin{enumerate}
\item
$[m,n-m]$.
\item
$[m-1,n-m+1]$, except when $\mu=(k,k)$.
\item
$[m-1,n-m,1]$.
\item
At least 2 different constituents of length~4, except when $l=2$, where
the product only has one constituent $[m-3,n-m+1,1^2]=[n-k-1,k-1,1^2]$ of length~4.
\item For $l=2$,  $[\mu]\cdot[\nu]$ always has the
    constituent $[m-2,n-m+1,1]=[n-k,k-1,1]$;
    for $k>3$ it also contains $[m-2,n-m,2]=[n-k,k-2,2]$, and
    for $k=3$, $n\geq 7$, it contains $[n-4,2^2]$.
\item When $\mu=(k,k)$ and $l>3$,   $[\mu]\cdot[\nu]$ has
at least 3 different constituents of length~4.
\item When $\mu=(k,k)$ and $l=3$,   $[\mu]\cdot[\nu]$ has
 a constituent $[k+1,k-1]$.
\item
For $\mu=(3^2)$, $\nu=(4,2)$ we have
$$[3^2]\cdot [4,2]=[5, 1] + [4, 1^2] + [3^2] + [3, 2, 1] + [2^2, 1^1]\:.$$
\end{enumerate}
Except for the product in (8),
the product $[\mu]\cdot [\nu]$ has at least five
(almost) extreme components.
\end{proposition}

\proof
Constituent (1) comes from Corollary~\ref{CDvir2}(i);
the constituents in (2) and (3) are obtained by applying Lemma~\ref{BK4.6}.
The constituents in (4) and (6)
are obtained using Corollary~\ref{CDvir2'}(ii) and
Theorem~\ref{Th1C}
and Theorem~\ref{Th2C}, respectively.
For (5), we apply Lemma~\ref{BK4.6} to $[\mu']\cdot[\nu]$, and then obtain
after conjugation constituents of  length~3 in $[\mu]\cdot[\nu]$
as given.
The constituent for (7) may be obtained with the help of \cite{V}
or from \cite{BWZ}.
Assertion (8)
{can easily be computed directly.

The final assertion follows by collecting in each case
suitable constituents found above.
The only critical situations to be discussed are the products
$[k,k]\cdot [n-3,3]$, where $k>3$.
Here, $m=k+3$ and we have found so far the constituents
$[k+3,k-3], [k+2,k-3,1]$ of (almost) maximal width and
two components of maximal length~4.
We just need to find a further component of almost maximal length~3.
Using Lemma~\ref{BK4.6},
we easily check that $[k,k]\cdot [k-1,k-1,2]$ has a constituent
$[2k-3,3]=[n-3,3]$, hence $[k-1,k-1,2]$ is a further constituent
of almost maximal length in $[k,k]\cdot [n-3,3]$.
}
\qed

\begin{cor}
Let $\mu, \nu \vdash n$ be proper 2-part partitions.
Then
$c=c([\mu]\cdot[\nu])\leq 4$
if and only if one of the following holds:
\begin{enumerate}[(1)]
\item $c=1$, when $n=2$, and the product is irreducible.
\item $c=2$, when $n=2k$ is even, and the product is  $[n-1,1]\cdot [k,k]$.
\item $c=3$, when $n=3$, $\mu=\nu=(2,1)$, or $n=4$, $\mu=\nu=(2^2)$.
\item $c=4$, when $n\geq 4$, $\mu=\nu=(n-1,1)$, or $n=2k+1$, $\{\mu,\nu\}=\{(n-1,1),(k+1,k)\}$,
or $n=6$, $\mu=\nu=(3^2)$.
\end{enumerate}
\end{cor}

\medskip

A {\em hook} partition is of the form $(a,1^b)$;
we call this partition a {\em proper hook} when
both $a>1$ and $b>0$ hold.

\begin{proposition}\label{hooks}
Let $\mu, \nu \vdash n$ be different proper hooks,
say $\mu=(n-k,1^k)$, $\nu=(n-l,1^l)$, with $1<l<k$
and $n-k>k$, $n-l>l$.
Set $\ga=\mu\cap\nu$; this is a partition of $m=n-k+l$.
Then we have the following constituents in $[\mu]\cdot[\nu]$:
\begin{enumerate}
\item
$[m,1^{n-m}]$.
\item
$[m-1,2,1^{n-m-1}]$.
\item
$[m-1,1^{n-m+1}]$.
\item
$[n-k-l,1^{k+l}]$.
\item
$[n-k-l+1,1^{k+l-1}]$.
\item
$[n-k-l,2,1^{k+l-2}]$.
\end{enumerate}
In all cases we have at least six (almost) extreme components.
\end{proposition}

\proof
The maximal width constituent in (1)
comes from Corollary~\ref{CDvir2}(i);
the almost maximal width constituents in (2) and (3)
are obtained using Lemma~\ref{BK4.6}.
The (almost) maximal length constituents in (4), (5) and (6)
come from conjugating the (almost) maximal width constituents in
$[\mu]\cdot [\nu']$.
Considering the length of the partitions in (1)-(6),
one easily sees that they are all different.
\qed
\medskip

\begin{cor}
Let $\mu, \nu \vdash n$ be hooks.
Then
$c=c([\mu]\cdot[\nu])\leq 4$
if and only if one of the following holds:
\begin{enumerate}[(1)]
\item $c=1$, when one of the partitions is a trivial hook $(n)$ or $(1^n)$.
\item $c=3$, when $\mu=\nu=(2,1)$.
\item $c=4$, when $n\geq 4$, and
$\mu,\nu \in \{(n-1,1),(2,1^{n-2})\}$.
\end{enumerate}
\end{cor}

\section{Key results on (almost) extreme components}\label{sec:tech}

From now on we fix the following notation:
\\
\centerline{$\mu, \nu \vdash n$\:,
$\mu \cap \nu = \gamma \vdash m$\:,
$d=n-m$\:.}

\smallskip

Our aim is to obtain information on (almost) extreme components;
we may (and will) focus on the components of (almost) maximal width
{since the corresponding results for components of (almost) maximal
length then follow by multiplying with the sign character}.
In the final section, we will see that these
results are strong enough
to prove the conjectured classification of the products with few components.

Because the special products considered
in the previous sections are sufficiently well understood (in particular
concerning this classification conjecture),
we may put these products aside and assume the following
properties of $\mu,\nu$, referred to as {\em Hypothesis $(*)$}:
\begin{center}
\begin{minipage}{12cm}
{
(1) $\mu, \nu \not\in \{ (n), (n-1,1), (1^n), (2,1^{n-2})\}$.
\\
(2) $\mu\neq \nu$, $\mu\neq \nu'$.
\\
(3)
$\mu,\nu$ are not both 2-line partitions or both hooks.
}
\end{minipage}
\end{center}
Here, by a {\em 2-line} partition we mean a partition which has
at most two rows or at most two
columns{, i.e., it is a 2-part partition or conjugate
to a 2-part partition}.
Note that we cannot have $\gamma=(m)$, since otherwise one of
$\mu$, $\nu$ is $(n)$. Also, the assumptions on $\mu,\nu$
imply that $m\geq 4$ and $m<n$.

\begin{lemma}\label{lem:4.6'}
Assume Hypothesis $(*)$.
Let one of the following situations be given:
\\
(i) $\nu / \gamma$ is a row.
\\
(ii) $[\mu / \gamma ]$ and $[\nu / \gamma]$ are irreducible and
$c([\mu / \gamma ] \cdot [\nu / \gamma])=2$.
\\
Let $\theta=(\theta_1,\theta_2,\dots)\vdash d+1$ be such that
$[\theta]$ appears as a constituent in
\begin{equation}\label{Expr}
\sum_{A\ \ga\mathrm{-node}}[\mu / \gamma_A] \cdot [\nu / \gamma_A] \:.
\end{equation}
Then in both situations above, we have $\theta_1\leq m-1$.
\end{lemma}

\proof
First note that  since $[\theta]$
appears in $[\mu/\ga_A]\cdot [\nu/\ga_A]$ for some removable $\ga$-node $A$,
it is a constituent of $[\rho]\cdot [\tau]$
for some constituents $[\rho]$ of  $[\mu/\ga_A]$ and $[\tau]$
of  $[\nu/\ga_A]$.
Then $\rho \subset \mu$, $\tau\subset \nu$ and hence
$\theta_1 \leq |\rho\cap\tau| \leq |\mu\cap\nu|=m$.
{Recall that $(*)$ implies $m\geq 4$;
hence if $d=1$, then $\theta_1 \leq d+1 = 2 \leq m-1$.
Thus we may assume that $d>1$.}

Now assume that $\theta_1 = m$.
Since $\rho\cap\tau \subseteq \mu\cap\nu$, this implies
$\rho\cap \tau = \ga$.

In Case (i),
$\nu / \ga_A$ is a union of a row and a node,
so $[\nu / \ga_A] \subseteq [d+1]+[d,1]$ (where the inclusion
here means a subcharacter).
As $\ga \subseteq \tau$ and $\ga \neq (m)$, we get $\tau=(d,1)$
and then $\ga=(m-1,1)$.
Since $\nu/\ga$ is a row of size $d\geq 2$, and $\nu \neq (n-1,1)$,
we must have $\nu=(m-1,d+1)$.
But then $|\theta|=d+1\leq m-1$, a contradiction.

In Case (ii), we are in the situation of Corollary~\ref{Cor:skewdia2C}(3).
Then one of the skew diagrams, say $\mu/\ga$,
is a nontrivial rectangle $(r^s)$, and the other is $(d-1,1)$
up to rotation and conjugation.
Then $[\rho]$ is one of $[r+1,r^{s-1}]$, $[r^s,1]$,
and $[\tau]$ is one of $[d,1]$, $[d-1,2]$, $[d-1,1^2]$ or their conjugates.
Assuming that $\ga \subset \tau$ is a hook, the conditions on the shapes of
$\mu/\ga$ and $\nu/\ga$ easily give a contradiction.
Hence $\ga$ is $(m-2,2)$ or its conjugate.
First let $m>4$; conjugating we may assume that $\ga=(m-2,2)$.
Then $\mu=(m-2,2^{s+1})$ and $\nu=((m-1)^2)$.
Thus $|\theta|=d+1=m-1$, again a contradiction.
If $m=4$, then up to conjugation we have $\mu=(2^{s+2})$ and $\nu=(n-3,3)$,
i.e., both partitions are 2-line partitions, a case we have excluded
above.
\qed

\medskip

We use Lemma~\ref{lem:prepalmext} and
Lemma~\ref{lem:4.6'} to obtain the following result
on non-special products that may have few components
of maximal width.
We know that the second case in Lemma~\ref{lem:4.6'}
only occurs when one of the skew characters
corresponds to a nontrivial
rectangle and the other one is $[d-1,1]$ or $[2,1^{d-2}]$;
conjugating, if necessary, we may assume that we are in
the first situation.
The lemma now gives information on components
of almost maximal width.

\begin{lemma}\label{ThDvir2'}
Assume Hypothesis $(*)$.
Let  $\hat\la = (\la_2, \la_3 , \ldots)\vdash d+1$
and $\la=(m-1,\hat\la)$.
\begin{enumerate}[(i)]
\item
Assume that $\nu / \ga$ is a row. \\
If $[\hat\la]$ appears with positive coefficient in the virtual character
$$
\chi=\sum_{A\ \ga\mathrm{-node}}[\mu / \gamma_A] \cdot [\nu / \gamma_A]
- [\mu/ \ga]\uparrow^{S_{d+1}} \:,$$
then $[\la]$ appears in $[\mu]\cdot[\nu]$, more precisely
$$ d(\mu,\nu;\la)= \langle \chi, [\hat\la] \rangle \; .$$
\item
Assume
$[\mu/\ga]=[\al]$ with $\al=(a^b)$, $a,b>1$,
and $[\nu/\ga]=[d-1,1]$. Set \;
$$\oal=(a+1,a^{b-2},a-1), \ual=(a^{b-1},a-1,1), \uoal=(a+1,a^{b-2},a-1,1)$$
and let $B_0$, $B_1$ denote the top and bottom
addable nodes for~$\al$.

If $[\hat\la]$ appears with positive coefficient
in the virtual character
$$\chi= \sum_{A} [\mu/\ga_A]\cdot [\nu/\ga_A]- ([\oal]+[\ual])\uparrow^{S_{d+1}}
$$
then $[\la]$ appears in $[\mu]\cdot[\nu]$, more precisely
$$ d(\mu,\nu;\la)= \langle \chi, [\hat\la] \rangle \; .$$
\end{enumerate}
\end{lemma}

\proof
By  Lemma~\ref{lem:4.6'}, in both cases
$\la_2 \leq m-1$, so $\la$ is a partition
and we can then use Lemma~\ref{lem:prepalmext}.
Hence, as $\la_1=m-1$ we obtain from Lemma~\ref{lem:prepalmext}:
$$
d(\mu,\nu;\la) =
\sum_{A \, \ga\text{-node}}
 \langle [\mu / \ga_A] \cdot [\nu/\ga_A] , [\hat{\la}]\rangle
  - \langle ([\mu/\ga]\cdot [\nu/\ga])\uparrow^{S_{d+1}},
  [\hat\la]\rangle\:.
$$
When $\nu/\ga$ is a row, $[\mu/\ga]\cdot [\nu/\ga]=[\mu/\ga]$,
and we have the statement in~(i).
In Case (ii),
$[\mu/\ga]\cdot [\nu/\ga]= [\al]\cdot [d-1,1]
=[\ual]+[\oal]$.
\qed
\medskip

We now want to get information on constituents in the product
of almost maximal width  when
$[\mu/\ga]\cdot [\nu/\ga]$ is homogeneous, i.e., there
is only one component of maximal width.
Based on finding constituents in the virtual character $\chi$
defined above, we have the following crucial result which will be applied
in the final section.

\begin{proposition}
\label{Lemma4.7'}\label{prop:4.7'}
Assume Hypothesis $(*)$.
Assume $[\mu / \gamma] = [\alpha]$ is irreducible and $\nu / \gamma$ is a row.
If there exists a removable $\gamma$-node $A_0$ disconnected
from $\nu /\gamma$
then
$$\chi=
\sum_{A\  \ga\mathrm{-node}}[\mu / \gamma_A] \cdot [\nu / \gamma_A]
- [\mu/ \ga]\uparrow^{S_{d+1}}\:$$
is a character and
one of the following holds.
\begin{enumerate}
\item $c(\chi)\geq 4$.
\item $c(\chi)=3$  and one of the following holds:
\begin{enumerate}
\item $d=2$, and we are not in one of the cases in (3).
\item $d=2k$ for some $k\in \N$, $k>1$ and
we have one of the following:
\\
$\mu=((a+k)^2)$, $\nu=(a^2,2k)$ for some $a>d$, \\
$\chi=[k+2,k-1]+[k+1,k]+[k+1,k-1,1]$,\\
or $\mu=(2^{a+k})$, $\nu=(2k+2,2^{a-1})$ for some $a>1$,\\
$\chi=[2^k,1]+[3,2^{k-2},1^2]+[2^{k-1},1^3]$,
\\
or  $\mu=((k+1)^3)$, $\nu=(3k+1,2)$, \\
$\chi=[k+1,k]+[k,k,1]+[k+1,k-1,1]$,
\\
or $\mu=(2^{k+2})$, $\nu=(2k+2,1^2)$, \\
$\chi=[3,2^{k-1}]+[2^k,1]+[2^{k-1},1^3]$.
\item $d>2$ and we have one of the following: \\
$\mu=((2d)^{a+1})$, $\nu=((2d)^a,d^2)$ for some $a\in \N$, \\
$\chi=2[d,1]+[d-1,1^2]+[d-1,2]$,
\\
or $\mu=((d+a+1)^d)$, $\nu=((d+a)^d,d)$ for some $a\in \N$, \\ $\chi=[2,1^{d-1}]+[2^2,1^{d-3}]+[3,1^{d-2}]$,\\
or $\mu=(d^{a+2})$, $\nu=(2d,d^a)$ for some $a\in \N$, \\
$\chi= [d,1]+[d-1,2]+[d-1,1^2]$.
\end{enumerate}
\item
$c(\chi)=2$
and one of the following holds:
\begin{enumerate}
\item
$d=1$.
\item $d=2$ and {we have one of the following:}\\
$\mu=(4^{a+1})$, $\nu=(4^a,2^2)$,
 $\chi=2[2,1]+[1^3]$,  \\
or
$\mu=(2^{a+2})$, $\nu=(4,2^a)$,
 $\chi=[2,1]+[1^3]$,
\\
or
$\mu=((a+3)^2)$, $\nu=((a+2)^2,2)$ for some $a\in \N$,
$\chi=[2,1]+[3]$.
\end{enumerate}
\end{enumerate}
Furthermore,
any constituent $[\theta]$ of $\chi$
gives a constituent $[m-1,\theta]$ in $[\mu]\cdot [\nu]$.
\end{proposition}
\smallskip

{\bf Remark.}
The case $c(\chi)=2$ may also occur when {$d>1$, }
$\chi=[d+1]+[d,1]$ and $\mu,\nu$ are both
hooks or both 2-part partitions.
But we had explicitly assumed in $(*)$ that
$\mu$, $\nu$ are not both hooks or both 2-part partitions.

\proof
We have already proved in Lemma~\ref{ThDvir2'} that every
constituent appearing with positive coefficient in $\chi$
gives a constituent in $[\mu]\cdot[\nu]$.

Let $A_0$ be a removable $\ga$-node, disconnected from $\nu/\ga$.
By assumption, we have
$$[\nu/\ga]=[d] \; , \;
[\nu/\ga_{A_0}]=[d+1]+[d,1]\:.$$

{\bf Case 1.} $A_0$ is disconnected from $\mu/\ga$.

Then
$$[\mu/\ga_{A_0}] = [\mu/\ga]\uparrow^{S_{d+1}}= [\al]\uparrow^{S_{d+1}}\:.$$
Consider
$$\chi_0=[\mu/\ga_{A_0}]\cdot [\nu/\ga_{A_0}] -  [\mu/\ga]\uparrow^{S_{d+1}}
= [\mu/\ga_{A_0}] \cdot [d,1]
= \sum_{B\ \al\text{-node}} [\al^B]\cdot [d,1]\:.$$
We may already note at this point that $\chi$ is then a character.
{In fact, this character is not homogeneous.
For $d=1$, it has exactly two components, and we are in case (3) above.
For $d>1$,}
one of the partitions $\al^B$ is not a rectangle,
and $\chi_0$ has more than two components, by Corollary~\ref{Cor:Nat}.
It has three components exactly if $d=2$; in this case, both
for $\al=(2)$ and $\al=(1^2)$ we get
$$\chi_0=[3]+2[2,1]+[1^3]\:.$$
Thus,  $\chi_0$ and hence
$$\chi=\chi_0+\sum_{A\  \ga\mathrm{-node} \atop A\neq A_0}[\mu / \gamma_A] \cdot [\nu / \gamma_A]$$
has at least four components, when $d>2$, and three components
when $d=2$.\\
Thus we are in situation (1) of the proposition for $d>2$ and in
situation (2) for $d=2$.
\smallskip

Having dealt with Case 1, we may now assume that we are in the
following situation:
\smallskip

{\bf Case 2.} Every removable $\ga$-node  disconnected from $\nu/\ga$
is connected to $\mu/\ga$.

Then $[\mu/\ga_{A_0}]$  has a constituent $[\al^{B_1}]$ for some addable
node $B_1$ of $\al$ by Lemma~\ref{lem:4.5'}.
Thus
$[\mu/\ga_{A_0}]\cdot [\nu/\ga_{A_0}]$
contains $[\al^{B_1}]\cdot ([d+1]+[d,1])$,
and thus it contains $[\al]\uparrow^{S_{d+1}}$.
Hence $\chi_0$ (as defined above) is a character and
we get constituents in the character $\chi$ from the character
$$\chi'=\sum_{A\  \ga\mathrm{-node} \atop A\neq A_0}[\mu / \gamma_A]
\cdot [\nu / \gamma_A]\:.$$
\smallskip

{\bf Case 2.1.}
Assume that  there is a further removable
$\ga$-node $A_1\neq A_0$ that is disconnected from
$\nu/\ga$.

Then as above
$\chi_1=[\mu/\ga_{A_1}]\cdot [\nu/\ga_{A_1}]$ (and hence $\chi$)
contains   $[\al^{B_2}]\cdot ([d+1]+[d,1])$, for some addable node
$B_2\neq B_1$,
and hence $[\al]\uparrow^{S_{d+1}}=\sum_B [\al^B]$.
This latter character is never homogeneous;
it has two components exactly when $\al$ is a rectangle
and three components exactly when $\al$ is a {\em fat hook}, i.e.,
it has exactly two different part sizes.
Otherwise we already get four components and thus $c(\chi)\geq 4$.

Now when $\al$ is a nontrivial rectangle, then $\al^{B_2}$ is not a rectangle,
and thus  $[\al^{B_2}]\cdot ([d+1]+[d,1])$ (and hence $\chi$)
has at least four components
by Corollary~\ref{Cor:Nat} (note that $|\al|\geq 4$).

When $\al$ is a trivial rectangle {and $d>1$},
we may still assume that $\al^{B_2}$ is not a rectangle, by interchanging $A_0$ and $A_1$, if necessary.
Then $[\al^{B_2}]\cdot ([d+1]+[d,1])$ (and hence $\chi$)
has again at least four components except if $d=2$,
when it has three components and then also $c(\chi)=3$.
{When $d=1$, i.e., $\al=(1)$, $\chi$ has exactly two components.}

\smallskip

{\bf Case 2.2.}
There is no removable $\ga$-node $\neq A_0$ disconnected from $\nu/\ga$, but
there is a removable $\ga$-node $A_1$ connected to $\nu/\ga$ that is disconnected
from $\mu/\ga$.

In this situation
$[\nu/\ga_{A_1}]$ is one of $[d+1]$ or $[d,1]$, and
$[\mu/\ga_{A_1}]= [\al]\uparrow^{S_{d+1}}$.

Clearly, $\chi_1=[\mu/\ga_{A_1}]\cdot [\nu/\ga_{A_1}]$ is not homogeneous.
{When $d=1$, we have $c(\chi)=2$, as required, so we may assume $d>1$.}
If $c(\chi_1)\geq 4$, we are done.
Now $c(\chi_1)=3$ if and only if $d=2$ and $[\nu/\ga_{A_1}]=[2,1]$,
or $\al$ is a fat hook
and $[\nu/\ga_{A_1}]=[d+1]$; in the first case, we are again done.
If $\mu/\ga=\al$ is a  fat hook,
$[\mu/\ga_{A_0}]$ has a further constituent $[\al^{B_2}]$, $B_2\neq B_1$,
and $c(\chi_0)\geq 4$.
If $\mu/\ga$ is a rotated fat hook,
$A_0$ must be its inner addable node.
As $\nu$ is not equal or conjugate to $(n-1,1)$,
$\chi_0$ has a constituent not appearing in
$\chi_1=[\al]\uparrow^{S_{d+1}}$,
and hence $\chi_0+\chi_1$ has at least four components.
Finally,
$c(\chi_1)=2$ if and only if
$[\nu/\ga_{A_1}]=[d+1]$ and $\al$ is a rectangle;
in this case, $[\mu/\ga_{A_0}]=[\al^{B_1}]$ and $\chi$ contains
$$\chi_0+\chi_1=[\al^{B_1}]\cdot ([d+1]+[d,1])\:.$$

If $\al^{B_1}$ is not a rectangle, then this has at least
four components and thus $c(\chi)\geq 4$,
except when $d=2$,  when we get $c(\chi)=3$.
If $\al^{B_1}$ is a rectangle, it must be a row or column.
If it is a row, then $\mu/\ga$ is a row, and the roles of $\mu$, $\nu$
can be interchanged. Thus we may assume that there is no
further removable $\ga$-node $\neq A_1$ disconnected from $\mu/\ga$.
But then we are in the situation where $\mu$, $\nu$
are both 2-part partitions:

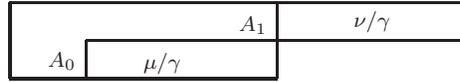
\begin{center}
\setlength{\unitlength}{0.01in}%
\begingroup\makeatletter\ifx\SetFigFont\undefined%
\gdef\SetFigFont#1#2#3#4#5{%
  \reset@font\fontsize{#1}{#2pt}%
  \fontfamily{#3}\fontseries{#4}\fontshape{#5}%
  \selectfont}%
\fi\endgroup%
\begin{picture}(220,100)(60,80)
\thicklines
\put( 60,180){\line( 1, 0){240}}
\put(100,160){\line( 1, 0){200}}
\put( 60,140){\line( 1, 0){140}}
\put( 60,180){\line( 0,-1){ 40}}
\put(100,160){\line( 0,-1){ 20}}
\put(200,180){\line( 0,-1) {40}}
\put(240,167){\makebox(0,0)[lb]{\smash{\SetFigFont{8}{7.2}{\rmdefault}{\mddefault}{\updefault}$\nu /\gamma$}}}
\put(300,180){\line( 0,-1){ 20}}
\put(130,145){\makebox(0,0)[lb]{\smash{\SetFigFont{8}{7.2}{\rmdefault}{\mddefault}{\updefault}$\mu /\gamma$}}}
\put(180,165){\makebox(0,0)[lb]{\smash{\SetFigFont{8}{7.2}{\rmdefault}{\mddefault}{\updefault}$A_1$}}}
\put(80,146){\makebox(0,0)[lb]{\smash{\SetFigFont{8}{7.2}{\rmdefault}{\mddefault}{\updefault}$A_0$}}}
\end{picture}
\end{center}
\vspace{-40pt}

If $\al^{B_1}$ is a column and there is a
further removable $\ga$-node $\neq A_1$ disconnected from $\mu/\ga$,
then we conjugate the partitions and use the previous arguments
to obtain at least four components in $\chi$, when $d>2$,
and three when $d=2$.
If there is no
further removable $\ga$-node $\neq A_1$ disconnected from $\mu/\ga$,
then $\mu$, $\nu$ are both hooks,
or we are in the following situation:
\\[5pt]

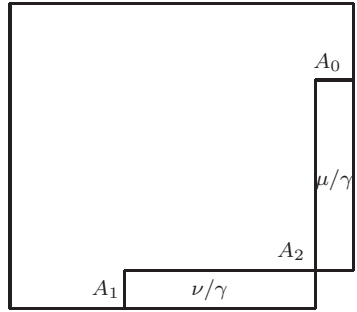
\begin{center}
\setlength{\unitlength}{0.01in}%
\begingroup\makeatletter\ifx\SetFigFont\undefined%
\gdef\SetFigFont#1#2#3#4#5{%
  \reset@font\fontsize{#1}{#2pt}%
  \fontfamily{#3}\fontseries{#4}\fontshape{#5}%
  \selectfont}%
\fi\endgroup%
\begin{picture}(220,80)(80,700)
\thicklines
\put( 80,800){\line( 1, 0){180}}
\put(240,760){\line( 1, 0){ 20}}
\put( 80,640){\line( 1, 0){160}}
\put( 140,660){\line( 1, 0){120}}
\put( 80,800){\line( 0,-1){160}}
\put(260,800){\line( 0,-1){ 140}}
\put(240,760){\line( 0,-1){ 120}}
\put(175,647){\makebox(0,0)[lb]{\smash{\SetFigFont{8}{7.2}{\rmdefault}{\mddefault}{\updefault}$\nu /\gamma$}}}
\put( 140,660){\line( 0, -1){ 20}}
\put(240,705){\makebox(0,0)[lb]{\smash{\SetFigFont{8}{7.2}{\rmdefault}{\mddefault}{\updefault}$\mu /\gamma$}}}
\put(123,647){\makebox(0,0)[lb]{\smash{\SetFigFont{8}{7.2}{\rmdefault}{\mddefault}{\updefault}$A_1$}}}
\put(220,667){\makebox(0,0)[lb]{\smash{\SetFigFont{8}{7.2}{\rmdefault}{\mddefault}{\updefault}$A_2$}}}
\put(239,766){\makebox(0,0)[lb]{\smash{\SetFigFont{8}{7.2}{\rmdefault}{\mddefault}{\updefault}$A_0$}}}
\end{picture}
\end{center}

\vspace{8ex}

In this latter case, we get a further contribution to $\chi$ from $A_2$:
$$\chi_2=[\mu/\ga_{A_2}]\cdot [\nu/\ga_{A_2}]
=[d,1]\cdot [d,1]
$$
and thus
$\chi=\chi_0+\chi_1+\chi_2$
has at least four components, except when $d=2$, where $c(\chi)=3$.
\smallskip

We now have to consider
\smallskip

{\bf Case 2.3.} All removable $\ga$-nodes $\neq A_0$ are connected with
$\nu/\ga$, and all removable $\ga$-nodes are connected with $\mu/\ga$.

Since $[\mu/\ga]$ is irreducible, $\mu/\ga$ is a partition or rotated partition.
\smallskip

{\bf Case 2.3.1} In the first case where $\mu/\ga$ is a partition
we have by Lemma~\ref{lem:4.5'}
$$
[\mu/\ga_{A_0}] = \sum_{B\neq B_0} [\al^B]
$$
where $B_0$ is the top or bottom addable node of $\al$.
Then
$$\begin{array}{rcl}
\chi_0&=&
\dstyle
[\mu/\ga_{A_0}]\cdot [\nu/\ga_{A_0}] - [\mu/\ga]\uparrow^{S_{d+1}}
= \left(\sum_{B\neq B_0} [\al^{B}]\right)([d+1]+[d,1])-\sum_{B} [\al^{B}]\\[8pt]
&=&
\dstyle [d,1] \cdot \left(\sum_{B\neq B_0} [\al^{B}]\right) - [\al^{B_0}] \:.
\end{array}
$$
As $d\geq 1$, $\al$ has a further addable node $B_1\neq B_0$.
If $\al$ is not a row or column, we may choose $B_1$ such that
$\al^{B_1}$ is not a rectangle.
Let $r$ be the number of removable nodes of $\al^{B_1}$.
Then by Lemma~\ref{LNat} we obtain
$$
\chi_0=
(r-1)[\al^{B_1}] + ([d+1]+[d,1])(\sum_{B\neq B_0,B_1} [\al^{B}])
+\sum_{{C\neq B_1} \atop {D\neq C}} [{(\al^{B_1})_C}^D]\:.
$$
We note that $\chi$ is thus a character,
and we consider the contributions coming from
the character
$$\chi_0'=([d+1]+[d,1])(\sum_{B\neq B_0,B_1} [\al^{B}])\:.$$

If $\al$ is not a fat hook  or rectangle,
then $\al$ has two further addable nodes
$B_2,B_3\neq B_0,B_1$.
Then $\chi_0'$ contains
$[\al^{B_0}]+[\al^{B_1}]+[\al^{B_2}]+[\al^{B_3}]$,
and thus $c(\chi)\geq 4$.

If $\al$ is a fat hook,
then it has a further
addable node $B_2\neq B_0,B_1$.
Then $\chi_0'$ contains $[\al^{B_0}]+[\al^{B_1}]+[\al^{B_2}]$,
and
$$\chi_0''=\sum_{{C\neq B_1} \atop {D\neq C}} [{(\al^{B_1})_C}^D]$$
contributes a further constituent as $\al^{B_1}$ has
a removable node $C\neq B_1$, and there is then  a suitable $D\neq C$
with ${(\al^{B_1})_C}^D \neq \al^{B_i}$, for $i=0,1,2$.

Now assume that $\al$ is a nontrivial rectangle, with corner node $Z$;
this is then the only removable node $\neq B_1$ of $\al^{B_1}$.
Then
$$\chi_0=\sum_{D} [{(\al^{B_1})_Z}^D]$$
has at least four different constituents,
except if $\al$ is a 2-row rectangle and $B_1$ is the top node,
or $\al$ is a 2-column rectangle and $B_1$ is the bottom node.
In these exceptional cases, if $B_1'$ is the top or bottom node
of $\al^{B_1}$, respectively, then
$\chi_0= [\al^{B_1}]+[(\al^{B_1})_Z^{B_0}]+[(\al^{B_1})_Z^{B_1'}]$
has exactly three components.

If there exists a removable $\ga$-node $A_1\neq A_0$, then this
must be connected to both $\mu/\ga$ and~$\nu/\ga$, and we have
$[\mu/\ga_{A_1}]=[\al^{B_0}]$ and $[\nu/\ga_{A_1}]$ is
$[d+1]$ or $[d,1]$.
In both cases, we find as a fourth new
constituent $[\al^{B_0}]$ in $\chi$.

Now it remains to consider the situation when there is no such $A_1$,
i.e., $\ga$ is a rectangle, and we are in one of the exceptional cases
where $d=2k$ for some $k>1$ and $\al$ is a 2-line rectangle.
Because of the additional condition on $B_1$,
we then have one of the following situations:

(i) $\ga=(2^a)$, $a\geq 2$, $\mu=(2^{a+k})$, $\nu=(2k+2,2^{a-1})$.

(ii) $\ga=(a^2)$, $a>d$, $\mu=((a+k)^2)$, $\nu=(a^2,2k)$.

Here, $c(\chi)=3$, and these situations appear
in part (2)(b) of the proposition.
\smallskip

Finally, we have  to deal with the case where $\al$ is a row or column.

First let $\al$ be a row.
Assume that there is a removable $\ga$-node $A_1\neq A_0$.
{If $d=1$, this leads to the contradiction that $\mu,\nu$
are both 2-line partitions. Hence we may now assume $d>1$.}
As $\mu,\nu$ are not both 2-part partitions,
$\mu/\ga_{A_0}$ cannot be a row,
but $\mu/\ga_{A_1}$ is a row, and hence we obtain
$$
\begin{array}{rcl}
\chi&=&
[d,1]\cdot ([d+1]+[d,1])+[d+1]\cdot [d,1]
-[d+1]-[d,1]
\\[8pt]
&=& \left\{
\begin{array}{ll}
2[d,1]+[d-1,2]+[d-1,1^2] & \text{if } d>2\\
2[2,1]+[1^3] & \text{if } d=2
\end{array}
\right. \:.
\end{array}
$$
We have here
$\mu=((2d)^{a+1})$, $\nu=((2d)^{a},d^2)$
for some $a\in \N$, cases described in the proposition.

If there is no such $\ga$-node $A_1$, then $\ga$ is  a rectangle.
{Because of Hypothesis $(*)$, we must then have $d>1$.}
Since  $\mu,\nu$ are not both 2-part partitions,
we must then have  $\mu=(d^{a+2})$, $\nu=(2d,d^a)$ for some $a\in \N$,
and we have
$$
\begin{array}{rcl}
\chi&=&
[d,1]\cdot ([d+1]+[d,1])-[d+1]-[d,1]
\\[8pt]
&=& \left\{
\begin{array}{ll}
[d,1]+[d-1,2]+[d-1,1^2] & \text{if } d>2\\
\-[2,1]+[1^3] & \text{if } d=2
\end{array}
\right. \:,
\end{array}
$$
cases appearing in the proposition.

Now let $\al$ be a column; we may assume $d>1$ since for $d=1$  we have that $\alpha$  is a row.
Assume that there is a removable $\ga$-node $A_1\neq A_0$.
If $\mu/\ga_{A_0}$ is a column,
we obtain
$$
\begin{array}{rcl}
\chi&=&
[1^{d+1}]\cdot ([d+1]+[d,1])+[d,1]\cdot [2,1^{d-1}]
-[1^{d+1}]-[2,1^{d-1}]
\\[8pt]
&=& [d,1]\cdot [2,1^{d-1}]\\[8pt]
&=& \left\{
\begin{array}{ll}
[1^{d+1}]+[2,1^{d-1}]+[2^2,1^{d-3}]+[3,1^{d-2}] & \text{if } d>2\\
\-
[3]+[2,1]+[1^3] & \text{if } d=2
\end{array}
\right. \:.
\end{array}
$$
This fits with the cases (1) and (2) in the proposition.

If $\mu/\ga_{A_0}$ is not a column,
then $\mu=((a+1)^{d+1})$, $\nu=(a+d+1,a^d)$ and
we obtain
$$
\begin{array}{rcl}
\chi&=&
[2,1^{d-1}]\cdot ([d+1]+[d,1])+ [1^{d+1}]
-[1^{d+1}]-[2,1^{d-1}]
\\[8pt]
&=& \left\{
\begin{array}{ll}
[1^{d+1}]+[2,1^{d-1}]+[2^2,1^{d-3}]+[3,1^{d-2}] & \text{if } d>2\\
\-[3]+[2,1]+[1^3] & \text{if } d=2
\end{array}
\right. \:.
\end{array}
$$
Again, this is in accordance with (1) and (2) of
the proposition.

When there is no removable $\ga$-node $A_1\neq A_0$,
then, since $\mu,\nu$ are not both hooks, we have
$\mu=((d+a+1)^d)$, $\nu=((d+a)^d,d)$, for some $a\in \N$,
and
$$
\begin{array}{rcl}
\chi&=& [2,1^{d-1}]\cdot ([d+1]+[d,1])
-[1^{d+1}]-[2,1^{d-1}]
\\[8pt]
&=& \left\{
\begin{array}{ll}
[2,1^{d-1}]+[2^2,1^{d-3}]+[3,1^{d-2}] & \text{if } d>2\\
\-[2,1]+[3] & \text{if } d=2
\end{array}
\right. \:.
\end{array}
$$
Again, these cases appear as exceptional situations in (2) and (3) of the
proposition.
\smallskip

{\bf Case 2.3.2.}
It remains to treat the case where $\mu/\ga$ is a rotated partition
which is not a partition.

Since only the removable $\ga$-node $A_0$ is disconnected from $\nu/\ga$,
$\mu/\ga$  can only be a  fat hook~$\al$,
and $A_0$ is the middle addable node $B_1$ (say) for~$\al$.
Then
$$\begin{array}{rcl}
\chi_0&=&[\mu/\ga_{A_0}]\cdot [\nu/\ga_{A_0}] -  [\mu/\ga]\uparrow^{S_{d+1}}
=[\al^{B_1}]\downarrow_{S_d}\uparrow^{S_{d+1}} -  [\al]\uparrow^{S_{d+1}}
\\[8pt]
&=&
\dstyle
\sum_{B\neq {B_1}\atop {D}} [{(\al^{B_1})_B}^D]\:,
\end{array}$$
where $B$ runs over the removable $\al^{B_1}$-nodes and
$D$ over the addable $(\al^{B_1})_B$-nodes;
in particular, we see here again that $\chi$ is a character.
There has to be exactly one further removable $\ga$-node $A_1$,
which corresponds to the top or bottom addable node $B_0$ or $B_2$ of $\al$,
respectively; in these two cases we obtain as the second
contribution to $\chi$:
$$
\chi_1=  [\mu/\ga_{A_1}]\cdot [\nu/\ga_{A_1}]
=\left\{
\begin{array}{l}
[\al^{B_0}]\cdot[d,1]=
\dstyle [\al^{B_1}]+[\al^{B_2}]+\sum_{B\neq {B_0}\atop {D}} [{(\al^{B_0})_B}^D]\\
\-[\al^{B_2}]
\end{array}
\right. \:.
$$
If $\al^{B_1}$ has three removable nodes $B_1,X,Y$,
then
$$\sum_{{D\neq X}} [{(\al^{B_1})_X}^D]+ \sum_{{D\neq Y}} [{(\al^{B_1})_Y}^D]$$
already gives at least four different constituents in $\chi_0$.

Now assume that $\al^{B_1}$ has only two removable nodes, $B_1$ and either
the top removable node $X$ or the bottom removable node $Y$ of $\al$.
Assume first that the top node $X$ is removable.
If $(\al^{B_1})_X$ has four addable nodes, then we have already
four different constituents in $\chi_0$.
If $(\al^{B_1})_X$ has three addable nodes,
then we have three different constituents from $\chi_0$
and a further fourth constituent
$[\al^{B_2}]$ from $\chi_1$ for $\chi$.
If $(\al^{B_1})_X$ has only two addable nodes,
then besides two constituents from $\chi_0$
we get at least two further constituents from $\chi_1$ when
$[\mu/\ga_{A_1}]=[\al^{B_0}]$.
When $[\mu/\ga_{A_1}]=[\al^{B_2}]=\chi_1$,
$\chi$ has only three components; in this situation we have
$d=2k$, $\mu=((k+1)^3)$, $\nu=(3k+1,2)$ for some $k\in \N$, $k>1$,
and $\chi=[k+1,k]+[k,k,1]+[k+1,k-1,1]$.

If the bottom node $Y$ is the second removable node of $\al^{B_1}$,
then we can argue analogously; this gives a further situation where
$\chi$ has three components, namely for
$d=2k$, $\mu=(2^{k+2})$, $\nu=(2k+2,1^2)$,
for some $k\in \N$, $k>1$; here $\chi=[3,2^{k-1}]+[2^k,1]+[2^{k-1},1^3]$.

Finally, we assume that $\al^{B_1}$ has only the removable node~$B_1$.
In this case $\chi_0=0$.
As by assumption $(*)$, $\nu\neq (n-1,1)$, the situation
$\chi_1=[\mu/\ga_{A_1}]=[\al^{B_2}]$ cannot occur.
Let $B_0'$ be the top addable node of $\al^{B_0}$.
Then $\chi=\chi_1$ has at least the four different
constituents $[\al^{B_i}]$, $i=0,1,2$, and
$[{(\al^{B_0})_Y}^{B_0'}]$.
\qed

\medskip

Now we want to deal with the second case in Lemma~\ref{lem:4.6'},
where we have a product $[\mu/\ga]\cdot [\nu/\ga]$ of two irreducible skew
characters with two components.
We know that this only occurs when one of the skew characters
corresponds to a nontrivial
rectangle, and the other one is
$[d-1,1]$ or $[2,1^{d-2}]$;
conjugating, if necessary, we may assume that we are in
the first situation.
The following result then provides lots of
components of almost maximal width in the product $[\mu]\cdot[\nu]$.

\begin{proposition}\label{prop:5.5}
Assume Hypothesis $(*)$.
Assume $[\mu/\ga]=[\al]$ with $\al=(a^b)$, $a,b>1$,
$[\nu/\ga]=[d-1,1]$ and
$$\chi=
\sum_{A\  \ga\mathrm{-node}}[\mu / \ga_A] \cdot [\nu / \ga_A]
- ([\mu/ \ga]\cdot [\nu / \ga])\uparrow^{S_{d+1}}\:.$$
Then $\chi$ is a character with $c(\chi)\geq 5$.

Furthermore,
any constituent $[\theta]$ of $\chi$
gives a constituent $[m-1,\theta]$ in $[\mu]\cdot [\nu]$.
\end{proposition}

\proof
Let $B_0$, $B_1$ denote the top and bottom
addable nodes for~$\al$; let $X$ be the removable $\al$-node.

Let $A$ be a removable $\ga$-node.
Then $[\mu / \ga_A]$ contains a constituent
$[\al^B]$, for $B=B_0$ or $B=B_1$  (for both
if $A$ is disconnected from $\mu/\ga$);
let $\bar B$ be the other addable node for~$\al$.
If $A$ is disconnected from $\nu / \ga$ or if
$\nu / \ga$ is a partition diagram,
then $[\nu / \ga_A]$ contains
$[d-1,1^2]+[d-1,2]$ or $[d-1,2]+[d,1]$.
If $\nu / \ga$ is a rotated partition and there is no
removable $\ga$-node disconnected from $\nu / \ga$, then we have
at least two removable $\ga$-nodes $A_0$ and $A_1$ connected to
$\nu/\ga$ giving us a contribution
$[\al^B]\cdot ([d-1,1^2]+[d-1,2])$
or  $[\al^B]\cdot ([d-1,2]+[d,1])$ to the sum in~$\chi$.

Thus we will now investigate the expressions
$$\chi'=[\al^B]\cdot ([d-1,1^2]+[d-1,2])- ([\al]\cdot[d-1,1])\uparrow^{S_{d+1}}$$
and
$$\chi'=[\al^B]\cdot ([d-1,2]+[d,1])- ([\al]\cdot[d-1,1])\uparrow^{S_{d+1}}\;,$$
respectively.
If then $\chi'$ is a character, so is $\chi$, and $c(\chi')\leq c(\chi)$.\\
In the following, instead of $\uparrow^{S_{d+1}}$ and similar
inductions one step up, we will just write $\uparrow$.

In the first case we use the relation
$[d-1,1^2]+[d-1,2]=[d-1,1]\uparrow -[d,1]$:
$$\begin{array}{rcl}
\chi'
&=&[\al^B]\cdot ([d-1,1^2]+[d-1,2])- ([\al]\cdot[d-1,1])\uparrow \\[8pt]
&=&
([(\al_X)^B]\cdot [d-1,1])\uparrow - [\al^B]\cdot[d,1]
\\[8pt]
&=& \dstyle
(r-1)[(\al_X)^B]\uparrow + \sum_{C}\sum_{D\neq C} [(((\al_X)^B)_C)^D]\uparrow
-[\al]\uparrow - [(\al_X)^B]\uparrow + [\al^B]
\end{array}
$$
where $r$ is the number of removable nodes of $(\al_X)^B$.
As $r\geq 2$, $\chi'$ contains the subcharacter
$$
\sum_{C\neq B}\sum_{ D\neq C} [({(\al_X)^B)}_C)^D]\uparrow
+\sum_{ D\neq B,X} [(\al_X)^D]\uparrow   + [\al^B]\:.
$$

Now $(\al_X)^B$ has at least one removable node $Y\neq B$, and
$\al_X$ has the addable node $\bar B\neq B,X$; hence the
character above  contains
$$
\sum_{ D\neq Y} [({(\al_X)^B)}_Y)^D]\uparrow
+[(\al_X)^{\bar B}]\uparrow   + [\al^B]\:,
$$
which has at least five components.
Thus in this case we have $c(\chi)\geq 5$.

Now we look at the second case. Since we have already dealt with the
previous case we may here assume that $B=B_1$.
When $\nu/\ga$ is a rotated partition, one of the two $\ga$-nodes
connected to $\nu/\ga$ is not connected to $\mu/\ga$, so that from
these nodes we get the subcharacter
$$
\begin{array}{rcl}
\chi'&=& [\al]\uparrow \cdot [d-1,2]+[\al^B]\cdot[d,1] -
([\al]\cdot[d-1,1])\uparrow
\\[8pt]
&=&
([\al]\cdot [d-2,2])\uparrow+[\al^B]\cdot[d,1]
\\[8pt]
\end{array}
$$
of $\chi$.
The second summand  has at least four constituents,
all of width $\leq a+1$,
and the first summand has one of width $a+3$, so
$\chi$ is a character with at least five
constituents in this case.
Now it only remains to consider the case where
$\mu=(a^{b+2})$ and $\nu=(a+d-1,a+1)$.
Since $\mu,\nu$ are not both 2-line partitions,
we have $a>2$.
We now want to show that the following is a subcharacter in $\chi$
with at least five components:
$$
\begin{array}{rcl}
\chi'&=&
[\al^B]\cdot ([d-1,2]+[d,1])- ([\al]\cdot[d-1,1])\uparrow
\\[8pt]
&=&
[\al^B]\cdot [d-1,2] +[\al^B]\cdot [d]\uparrow - [\al^B]
-([(\al_X)^B] + [(\al_X)^{\bar B}])\uparrow
\\[8pt]
&=&
[\al^B]\cdot [d-1,2] + [\al]\uparrow -  [\al^B] -[(\al_X)^{\bar B}]\uparrow
\\[8pt]
&=&
[\al^B]\cdot [d-1,2] + [\al^{\bar B}] -[(\al_X)^{\bar B}]\uparrow
\\[8pt]
&=&
[\al^B]\cdot [d-1,2]  -\sum_{C\neq X}[((\al_X)^{\bar B})^C]
\;.
\end{array}
$$
Note that the sum that is subtracted above has at most three terms,
namely $[a+2,a^{b-2},a-1]$ and $[a+1,a^{b-2},a-1,1]$,
and for $b \geq 3$ also $[(a+1)^2,a^{b-3},a-1]$.

We now investigate the product $[\al^B]\cdot [d-1,2]$.
As $B=B_1$,  $\al^B\cap (d-1,2) =(a,2) =\tau $ and thus
$$[\al^B/\tau]\cdot [(d-1,2)/\tau ] = [(a^{b-2},1)/(2)]
= [a^{b-2},a-2,1]+ [a^{b-2},a-1] \:,
$$
giving the components
$[a+2,a^{b-2},a-2,1]$, $[a+2,a^{b-2},a-1]$
in $[\al^B]\cdot [d-1,2]$. \\
We compute the terms of width $a+1$ in  $[\al^B]\cdot [d-1,2]$
to see that $\chi'$ is a character.
Since $[(d-1,2)/\tau ] $ is a row, we can use Lemma~\ref{ThDvir2'},
so we now compute the constituents of
$$
\begin{array}{rcl}
\psi &= & [(a^{b},1)/(a-1,2)] +  [(a^{b-1},1)/(1)]\cdot [d-2-a]\uparrow -
([(a^{b-2},a-2,1)]+ [a^{b-2},a-1])\uparrow \\[8pt]
&= & [(a^{b},1)/(a-1,2)] +  [(a^{b-1},1)/(1)]\downarrow\uparrow -
([(a^{b-2},a-2,1)]+ [a^{b-2},a-1])\uparrow \:.
\end{array}
$$
Now for the first term in $\psi$ we have for $a\geq 4$ (see \cite{CG})
$$[(a^{b},1)/(a-1,2)] =[(a^{b-1},a-2,1)/(a-1)]= [a^{b-2},a-2,1^2]+[a^{b-2},a-2,2]+[a^{b-2},a-1,1]$$
while for $a=3$ the second summand does not appear, i.e.,
$$ [(3^b,1)/(2^2)]=[(3^{b-1},1^2)/(2)]=[3^{b-2},1^3]+[3^{b-2},2,1]\:.$$
For the second term in $\psi$, we first get
$$[(a^{b-1},1)/(1)]
=[a^{b-2},a-1,1]+[a^{b-1}]\:.
$$
Now we notice that the restriction of the first summand
already contains the
two constituents in the third term subtracted in the expression for~$\psi$.
Then from the second and third term in $\psi$ together we obtain
the contribution
$$
([a^{b-3},(a-1)^2,1]+[a^{b-2},a-1])\uparrow \:,
$$
where the first constituent only appears for $b\geq 3$.

Hence $\psi$ is a character, and taking into regard the
contribution from the first term, it contains
for all $a\geq 3$, $b\geq 2$ the character
$$
\begin{array}{rcl}
\psi' & = & [a^{b-2},a-2,1^2]+2[a^{b-2},a-1,1]
+[a+1,a^{b-3},a-1]+[a^{b-2},a]\:.
\end{array}
$$
All these constituents in $\psi'$ give constituents of width $a+1$ in
$[\al^B]\cdot [d-1,2]$, and thus the subtracted terms in the expression
for $\chi'$ are all taken care of, i.e., $\chi'$ is a character and
it contains the character
$$\chi''=
[a+2,a^{b-2},a-2,1]+
[a+1,a^{b-2},a-2,1^2]+[a+1,a^{b-2},a-1,1]+[a+1,a^{b-2},a]\:.
$$
Furthermore, as $d=ab \geq b+2$ and hence $d-1\geq b+1$,
$\al^B\cap (2^2,1^{d-3}) =(2^2,1^{b-1}) =\rho $ and thus
$$[\al^B/\rho]\cdot [(2^2,1^{d-3})/\rho ] =
[(a-1)^{b-2},(a-2)^2]\cdot [1^{d-b-2}]
= [b^{a-2},b-2]
$$
producing the only component of maximal length $b+3$ in
$[\al^B]\cdot [d-1,2]$:
$[a^{b-2},(a-1)^2,1^3]$, then also appearing in~$\chi'$.
Thus we have proved that $\chi$ is a character and
$c(\chi)\geq c(\chi') \geq 5$.
\qed

\section{Proof of the classification theorems}\label{sec:proof}

We are now in a position to
prove {the classification stated in Theorem~\ref{the:extcomp} as well as
the conjectured classification of Kronecker products with
only three or four homogeneous components in Theorem~\ref{the:34C}}.

\medskip

{We start with the {\bf Proof of Theorem~\ref{the:extcomp}.}}

We keep the notation used in earlier sections, i.e., we
set
$$\ga=\mu\cap\nu \vdash m, d=n-m.$$

We recall Hypothesis $(*)$ from Section~\ref{sec:tech}:
\begin{center}
\begin{minipage}{12cm}
{
(1) $\mu, \nu \not\in \{ (n), (n-1,1), (1^n), (2,1^{n-2})\}$.
\\
(2) $\mu\neq \nu$, $\mu\neq \nu'$.
\\
(3)
$\mu,\nu$ are not both 2-line partitions or both hooks.
}
\end{minipage}
\end{center}
Since we know the result for the cases
excluded in~$(*)$ by Corollary~\ref{Cor:Nat}
and Propositions~\ref{prop:TensSqEx}, \ref{twopart}
and~\ref{hooks},
we may (and will) assume that Hypothesis $(*)$ is
satisfied for the given partitions~$\mu,\nu$.

Note that all the exceptional cases
on the classification list
are pairs of partitions put aside when assuming
Hypothesis $(*)$.
Thus we are now in the situation that we
want to find at least five (almost) extreme components
in the product $[\mu]\cdot [\nu]$.

\smallskip

We note at this point that Hypothesis $(*)$
above implies
$|\mu\cap \nu'| \geq 4$,
and thus a maximal length component of $[\mu]\cdot [\nu]$
is of length~$\geq 4$.

\subsection{One component of one extreme type}\label{subsec:one}

We consider first the situation
that for one of the two extreme types
(maximal width or maximal length) there is
only one component in the product.

Replacing, if necessary, one of the partitions by its conjugate,
we may assume that there is only one component of maximal width~$m$.
Then, by Corollary~\ref{CDvir2}(i),
we know that
$[\mu/\gamma]\cdot[\nu/\gamma]$
must be homogeneous.
By Corollary~\ref{Cor:skewdia1C},
both skew characters have to be irreducible,
and one of them is of degree~1.
Conjugating both partitions and renaming,
if necessary, we may then assume that
$[\mu/\gamma]$ is irreducible and that $\nu/\gamma$ is a row.
\smallskip

If there is a removable $\ga$-node $A_0$ disconnected from $\nu/\ga$,
we can use Proposition~\ref{Lemma4.7'} to obtain constituents of
almost maximal width $m-1$
in the product.
Let $$\chi=
\sum_{A\  \ga\mathrm{-node}}[\mu / \gamma_A] \cdot [\nu / \gamma_A]
- [\mu/ \ga]\uparrow^{S_{d+1}}$$
be as before.
By Proposition~\ref{Lemma4.7'},
we obtain at least four components of width $m-1$
in the product coming from constituents of $\chi$,
unless we are in one of the exceptional cases described explicitly in
Proposition~\ref{prop:4.7'}.
 We now go through these in detail.
\smallskip

Assume first that $c(\chi)=3$.

First we consider the case $d=2$, where
we have already found one of the constituents
$[n-2,2]$ or $[n-2,1^2]$,  and the constituents
$[n-3,3]$, $[n-3,2,1]$, $[n-3,1^3]$ in the product.
If $|\mu\cap\nu'|>4$, then we also have a constituent of maximal
length $>4$.
Thus we may now assume that  $|\mu\cap\nu'|=4$.
In this case $(3,1^2)$ and $(2^2)$ cannot be contained in~$\ga$,
hence $\ga=(2,1^{m-2})$ or $\ga=(m-1,1)$.
Hypothesis $(*)$ then gives a contradiction except for the
case
$\mu=(2^3)$, $\nu=(4,1^2)$ or the (doubly conjugate) case
$\mu=(3,1^3)$, $\nu=(3^2)$.
Since $[\mu/\mu\cap \nu']\cdot [\nu'/\mu\cap \nu']= [2]+[1^2]$,
we get here a second constituent $[2^2,1^2]$ of length~4
in~$[\mu]\cdot[\nu]$.
Hence we have found five (almost) extreme components in the product.

Next suppose $d=2k$ for some $k>1$.
If $\mu=((a+k)^2)$, $\nu=(a^2,2k)$ for some $a>d$,
then $\chi=[k+2,k-1]+[k+1,k]+[k+1,k-1,1]$,
and we have four components of (almost) maximal width, and
of length $\leq 4$.
Since $|\mu\cap\nu'|=6$, we also have
an extreme component of length~6.
Similarly, when
$\mu=((k+1)^3)$, $\nu=(3k+1,2)$ for some $k\in \N$,  $k>1$,
we have $\chi=[k+1,k]+[k,k,1]+[k+1,k-1,1]$,
and thus we have again four
components of (almost) maximal width and of length $\leq 4$,
and a further extreme
component of length $|\mu\cap\nu'|=5$ in the product.
If $\mu=(2^{a+k})$, $\nu=(2k+2,2^{a-1})$ for some $a>1$,
then $\chi=[2^k,1]+[3,2^{k-2},1^2]+[2^{k-1},1^3]$, and
we have four (almost) extreme components of length
$\leq k+3$ and
a fifth extreme component of length $|\mu\cap \nu'|\geq k+4$.
When $\mu=(2^{k+2})$, $\nu=(2k+2,1^2)$,
for some $k\in \N$, $k>1$,
we have $\chi=[3,2^{k-1}]+[2^k,1]+[2^{k-1},1^3]$,
hence we have three (almost) extreme components of length $\leq k+2$
and the extreme component
$[3,2^{k-1},1^3]$ of length~$k+3=|\mu\cap \nu'|$.
Now we have to look more closely
at the components of maximal width
in $[\mu]\cdot[\nu']$.
Since\\
\centerline{$
[\mu/(\mu\cap \nu')]\cdot [\nu/(\mu'\cap \nu)] =
[1^{k+1}]+[2,1^{k-1}]
$}\\[5pt]
by Corollary~\ref{CDvir2'}(ii), the product $[\mu]\cdot [\nu]$
has $[2^{k+1},1^2]$ and $[3,2^{k-1},1^3]$
as components of maximal length.
Hence we have found five (almost) extreme components in the product.

Now consider the cases for $d>2$ in Proposition~\ref{prop:4.7'}(2)(c).
When $\mu=((2d)^{a+1})$, $\nu=((2d)^a,d^2)$ for some $a\in \N$,
we have $\chi=2[d,1]+[d-1,1^2]+[d-1,2]$, hence there are
already four (almost) extreme components of length $\leq 4$,
and because $|\mu\cap\nu'|\geq 6$,
we also have am extreme component of length~$\geq 6$.
Similarly, if $\mu=(d^{a+2})$, $\nu=(2d,d^a)$ for some $a\in \N$,
we have $\chi=[d,1]+[d-1,2]+[d-1,1^{2}]$.
Thus we have four (almost) extreme components of length $\leq 4$, and
a further extreme component of length $|\mu\cap \nu'|\geq 6$.
When $\mu=((d+a+1)^d)$, $\nu=((d+a)^d,d)$ for some $a\in \N$,
we have $\chi=[2,1^{d-1}]+[2^2,1^{d-3}]+[3,1^{d-2}]$.
Thus we have four (almost) extreme components of length $\leq d+1$, and
a further extreme component of length $|\mu\cap \nu'|=d(d+1)> d+1$.

Next we consider the cases where $c(\chi)=2$.
We start with the cases for $d=2$.\\
Consider the case
$\mu=(4^{a+1})$, $\nu=(4^a,2^2)$
where $\chi=2[2,1]+[1^3]$; then
$[\mu]\cdot [\nu]$ has three (almost) extreme components
$[4a+2,2]$, $[4a+1,2,1]$
and $[4a+1,1^3]$ of length $\leq 4$.
For $a=1$, $\mu'$ and $\nu$ satisfy the assumptions
of Proposition~\ref{prop:4.7'}, hence
the product $[\mu']\cdot [\nu]$
has three components of (almost) maximal width,
giving three components of (almost) maximal
length $\geq 5$ in $[\mu]\cdot [\nu]$.
For $a=2$, $\nu=\nu'$, hence $[\mu]\cdot [\nu]=[\mu]\cdot [\nu']$
has three components of (almost) maximal length $\geq 4a+1=9$;
similarly, for $a=3$, $\mu=\mu'$ and thus
the product has three components of (almost) maximal length $\geq 4a+1=13$.
For $a=4$, the pair $\mu$, $\nu'$ satisfies the assumption
of Proposition~\ref{prop:4.7'} (with the row diagram
$\mu/(\mu\cap \nu')$); hence $[\mu]\cdot [\nu']$
has three components of (almost) maximal width $\geq |\mu\cap \nu'|-1=15$,
giving three components of (almost) maximal
length $\geq 15$ in $[\mu]\cdot [\nu]$.
Finally, for $a>4$, $[\mu]\cdot [\nu']$
has at least two components of
maximal width~16, hence  $[\mu]\cdot [\nu]$ has
at least two components of maximal length~16.

Now suppose
$\mu=(2^{a+2})$, $\nu=(4,2^a)$;
here $\chi=[2,1]+[1^3]$, so we have again three
(almost) extreme components $[2a+2,2]$, $[2a+1,2,1]$ and
$[2a+1,1^3]$  of length $\leq 4$ in the product.
As $\mu,\nu$ are not both 2-line partitions, we have $a>1$.
For $a=2$, we may use Proposition~\ref{prop:4.7'}
to get three components of (almost) maximal width $\geq 5$ in
$[\mu']\cdot [\nu]$, and hence
three components of (almost) maximal length $\geq 5$ in
$[\mu]\cdot [\nu]$.
For $a\geq 3$, $[\mu']\cdot [\nu]$
has at least two components of maximal width~6,
giving two components of maximal length~6
in $[\mu]\cdot [\nu]$.

Finally, consider the case
$\mu=((a+3)^2)$, $\nu=((a+2)^2,2)$ for some $a\in \N$;
here $\chi=[2,1]+[3]$, so we have three
components $[2a+4,2]$, $[2a+3,2,1]$, $[2a+3,3]$
of (almost) maximal width and of length $\leq 3$
in $[\mu]\cdot [\nu]$.
For $a=1$, $[\mu]\cdot [\nu']$ has $\geq 3$ components of
(almost) maximal width
$\geq 5$ by Proposition~\ref{prop:4.7'}, hence
$[\mu]\cdot [\nu]$ has $\geq 3$ components of (almost) maximal length
$\geq 5$.
For $a>1$,
$[\mu/(\mu\cap \nu')]\cdot [\nu'/(\mu\cap \nu')]=[a^2]^2$
has $\geq 3$ components by Proposition~\ref{LTensSq},
hence $[\mu]\cdot [\nu]$ has $\geq 3$ components of maximal length
$|\mu\cap \nu'|=6$.

It remains to discuss the case where $d=1$.
Here, we have the three (almost) maximal width constituents $[n-1,1]$,
$[n-2,1^2]$ and $[n-2,2]$, as well as an extreme component
of length~$|\mu \cap \nu'|\geq 4$.
If there is a second component of maximal length, we are done.
Thus it remains to discuss the case where
{
$[\mu/\mu\cap\nu']\cdot[\nu'/\mu\cap\nu']$
is homogeneous, and hence irreducible, by Corollary~\ref{Cor:skewdia1C}}.
Then one of the skew diagrams
$\mu/\mu\cap\nu',\nu'/\mu\cap\nu'$
is a row or column
and the other is a partition or rotated partition;
conjugating both partitions, if necessary, we have
a row $D$ and a (rotated) partition.
If there is a removable $\ga$-node $A_1$ disconnected from
the row, then we obtain at least two components of
$[\mu]\cdot[\nu']$ of almost maximal width, and hence
two components of $[\mu]\cdot[\nu]$ of almost maximal
length~$|\mu \cap \nu'|-1\geq 3$;
thus we have at least five (almost) extreme components in the
product.
Hence we may now assume that any removable $\ga$-node
is connected to the row $D$; then $D$ complements
its partition to a rectangle, say $(a^b)$, where
$a,b >2$.
Since $\mu,\nu$ differ only by moving one node, and there is
a removable $\ga$-node
disconnected from $\nu/\ga$, we can then only have
for the two partitions,
up to conjugation, $\mu=(a^a), \nu=(a+1,a^{a-2},a-1)$,
or $\mu=(a^{a+1})$ and  $\nu\in \{(a+1,a^{a-1},a-1),(a^a,a-1,1)\}$.
In the first case, we get by symmetry six (almost) extreme components.
In the other two cases, with Lemma~\ref{lem:prepalmext} we also
get at least two components of (almost) maximal length in $[\mu]\cdot[\nu]$.
So we are done in this case.

\medskip

The critical situation to be discussed now
is the one where all removable $\ga$-nodes are connected
to~$\nu/\ga$.
In this case, $\nu$ must be a rectangle,
since $\nu/\gamma$ is a row; since $\nu\neq (n)$,
$\ga$ must have a removable node $A_0$ such that
$[\nu/\ga_{A_0}]=[d,1]$.
{
When $d=1$, we may interchange the partitions $\mu,\nu$,
and thus may also assume that all removable $\ga$-nodes  are connected
to $\mu/\ga$. But then $\mu,\nu$ are both 2-part partitions,
contradicting Hypothesis $(*)$.
So we can now assume that $d>1$.
}

Let us first assume that $\mu / \ga$
is disconnected from~$A_0$.
Then
$$\label{XX}
\begin{array}{rcl}
\chi_0&=&
\dstyle
[\nu / \gamma_{A_0}] \cdot [\mu / \gamma_{A_0}]
- \sum_B [\alpha^B]
 =  \displaystyle
[d,1]   \cdot \Big(\sum_B [\alpha^B]\Big) - \sum_B [\alpha^B]\\[5pt]
& = & \displaystyle \sum_B \sum_{C} \sum_D [{(\alpha^B)_C}^D]
- 2 \sum_B [\alpha^B] \\[5pt]
& = & \displaystyle \sum_B (r_B-2) [\alpha^B] +
      \sum_B \sum_{C} \sum_{D \neq C} [{(\alpha^B)_C}^D]
\end{array}
$$
where $B$ runs over the addable nodes of $\al$,
$C$ runs over the removable nodes of $\al^B$ (for the
respective node $B$), $D$ runs over the addable nodes of
$(\al^B)_C$ and $r_B$ denotes the number of
removable nodes of $\al^B$.

We know  that $\al$
has at least two addable nodes, say
$B_0$ at the top and $B_1$ at the bottom.
\\
Assume first that $\al$ is not a row or column.
Then $r_{B_0}, r_{B_1} \geq 2$.
Let $X_1$ be the top and $X_0$ the bottom
removable node of $\al$ (we may have $X_0=X_1$);
then $X_i$ is also a removable node for $\al^{B_i}$,
$i=0,1$. Let $B_0'$ be the top addable
node  for $(\al^{B_0})_{X_0}$ and
$B_1'$ be the bottom addable
node  for $(\al^{B_1})_{X_1}$.
We then have at least the following
contribution to $\chi_0$:
$$
(r_{B_0}-1)   [\al^{B_0}] +  [{(\al^{B_0})_{X_0}}^{B_0'}]
+(r_{B_1}-1)   [\al^{B_1}]+[{(\al^{B_1})_{X_1}}^{B_1'}]\:.
$$
Thus we have found at least
four components in $\chi_0$ and hence
four components of almost maximal width~$m-1$ in
the product.

If $\al$ is a row, we may interchange $\mu$ and $\nu$, and then we are in the situation discussed in the first part of the proof.
When $\al$ is a column, we conjugate and then interchange
both partitions; again, this is dealt with
by the first part of the proof.

\smallskip

Thus now we treat the situation where $\mu/\ga$ is connected to $A_0$.
Then by Lemma~\ref{lem:4.5'}
$$[\mu / \ga_{A_0}] = \sum_{B \neq B_1} [\al^B]$$
where $B_1$ is the bottom addable node of $\al$.
Let $B_0$ be the top addable node of $\al$.
Then
$$
\begin{array}{rcl}\label{XXX}
\chi_0
&= & [\nu / \ga_{A_0}] \cdot [\mu / \gamma_{A_0}]
- \dstyle \sum_B [\al^B]
=
[d,1]  \cdot \sum_{B \neq B_1} [\al^B] - \sum_B [\al^B]
\\[10pt]
& = & \dstyle \sum_{B \neq B_1} \sum_{C} \sum_D [{(\al^B)_C}^D]
- \sum_{B \neq B_1} [\al^B] - \sum_B [\al^B]
\\[10pt]
& = & \dstyle \sum_{B \neq B_0,B_1} \sum_{C} \sum_D [{(\al^B)_C}^D]
+ \sum_{C \neq B_0} \sum_D [{(\al^{B_0})_C}^D]
- \sum_{B \neq B_1} [\al^B]
\end{array}
$$
where $B$ runs through the addable nodes of $\al$,
$C$ runs through the removable nodes of $\al^B$ (for the
respective node $B$) and $D$ runs through the addable nodes of
$(\al^B)_C$.

If $\al$ is not a rectangle,
there is a third addable node, say $B_2$,
not in the first row or column.
Taking this contribution into account,
$\chi_0$ is a character containing
$$
\chi_0'=
[\al^{B_0}]+ \sum_{C \neq B_0} \sum_{D\neq C} [{(\al^{B_0})_C}^D]
+ \sum_{C \neq B_2}\sum_D [{(\al^{B_2})_C}^D]
+ [\al^{B_1}]\:.
$$

If $\al^{B_2}$ is not a rectangle,
then the top or bottom removable
node of $\al$ will be  $\neq B_2$ and will
also be removable from $\al^{B_2}$;
let this $\al$-node be~$X$.
Depending on $X$ being at the top or bottom of~$\al$,
the node $Y=B_1$ or $Y=B_0$ will be addable for $(\al^{B_2})_X$.
Then $\chi_0'$ contains
$$
\chi_0''=
[\al^{B_0}]+  [\al^{B_2}]
+ [{(\al^{B_2})_X}^Y]
+ [\al^{B_1}]
$$
and hence we have at least four components in the product of almost maximal width~$m-1$.

We are now in the situation where $\al^{B_2}$ is a rectangle.
The bottom removable node $X_0$ of~$\al$ is also removable
from~$\al^{B_0}$. The top addable node $B_0'$ for  $\al^{B_0}$
is also addable for $(\al^{B_0})_{X_0}$.
If $\al \neq (2,1)$, then $(\al^{B_0})_{X_0}$
has a second addable node $Y=B_1$ or $Y=B_0''$ (in the second row).
Thus, in this situation $\chi_0'$ contains
$$
\chi_0''=
[\al^{B_0}]+  [{(\al^{B_0})_{X_0}}^{B_0'}]
+ [{(\al^{B_0})_{X_0}}^Y]
+ [\al^{B_1}]\:,
$$
giving us again four components of almost maximal width~$m-1$ in the product.

When $\al=(2,1)$, we have
$$\chi_0=
[\al^{B_0}]+  [{(\al^{B_0})_{X_0}}^{B_0'}]
+ [\al^{B_1}]=[3,1]+[4]+[2,1^2]\:,
$$
so that up to this point we have found
four (almost) extreme components of length $\leq 4$ in $[\mu]\cdot[\nu]$.

If there is a removable $\ga$-node $A_1\neq A_0$ connected to $\nu/\ga$,
then we also get a contribution to $\chi$ from
$$\chi_1=[\mu/\ga_{A_1}]\cdot [\nu/\ga_{A_1}]
=[\al]\uparrow^{S_4}=[3,1]+[2^2]+[2,1^2]\:,$$
and thus we have again four components of almost maximal width~$m-1$
in the product.

If there is no such $\ga$-node $A_1$,
then $\mu=(5,4)$, $\nu=(3^3)$,
and $|\mu\cap\nu'|=|\mu\cap \nu|=6$
yields a component of maximal length~$6$ in the product.

Next we have to consider the case when $\al$ is a rectangle;
let $B_0$ be the top and $B_1$ the bottom addable node of~$\al$.
Let $X$ be the corner node of $\al$.
Since $\mu, \nu$ are not both 2-part partitions,
$\al$ is not a row.
Then $\al^{B_0}$ also has the removable
node~$X$ and we have
$$\chi_0=
\sum_D  [{(\al^{B_0})_X}^D] - [\al^{B_0}]
= \sum_{D \neq X}  [{(\al^{B_0})_X}^D]\:,$$
which gives three components of almost maximal width~$m-1$,
except in the cases where $\al$ has only two rows or
only one column.
When $\al=(1^2)$, we
have $\chi_0=[3]$, and otherwise, when $\al\neq (1^2)$
has only two rows
or one column, $\chi_0$ has two constituents.

Now if $\ga$ has a further removable node $A_1$, then
we also get two components of almost maximal width $m-1$ from
$$\chi_1=[\mu/\ga_{A_1}]\cdot[\nu/\ga_{A_1}]
=\sum_B [\al^B]=[\al^{B_0}]+[\al^{B_1}]\:,$$
and these are different from the ones appearing in $\chi_0$.
Thus for $\al\neq (1^2)$, we have then found at least
four components of almost maximal width~$m-1$ in the product.
For $\al=(1^2)$ we have found at this stage
four (almost) extreme components of length $\leq 4$;
but here $|\mu\cap\nu'| \geq 6$, giving us
also a component of maximal length $\geq 6$ in the product.

Thus we are now in the situation where $\ga$ is a rectangle
and $\al$ is a rectangle, say $\al=(a^b)$,
and we need to find  further (almost) extreme components.
We already know that $\al$ is not a row, and it also
cannot be a column because this would
contradict $\mu\neq \nu'$, so $1<a,b<d$, $d\geq 4$.

First assume $b\geq 3$.
By the above, we already have four (almost) extreme components
in the product which are of length $\leq b+2$.
Here $\mu\cap \nu'=((b+1)^b)$, hence we also get
an extreme component of length $b(b+1) > b+2$.

It remains to consider the case $b=2$,
i.e., $\al=(a^2)$ and $\mu=((3a)^2)$, $\nu=((2a)^3)$,
where $a>1$.
By the considerations so far, we have found
the (almost) extreme constituents $[4a,a^2]$, $[4a-1,a+2,a-1]$
and $[4a-1,a+1,a-1,1]$
in the product $[\mu]\cdot[\nu]$
which are of length $\leq 4$.
Now $\mu\cap \nu'=(3^2)$ and
$$[\mu /(\mu\cap \nu')]\cdot [\nu'/(\mu\cap \nu')]
= [(3a-3)^2]\cdot [3^{2a-2}]\:.
$$
By Theorem~\ref{Th1C} and Theorem~\ref{Th2C},
this product has at least three components, hence
by Corollary~\ref{CDvir2'}(ii)   $[\mu]\cdot[\nu]$
also has at least three extreme components of length~6.

{
Thus at this stage we have proved our claim for the case
that there is only one component for one of the two
extreme types.
}
\medskip

\subsection{Two components of each extreme type}\label{subsec:two}

{Towards the proof of Theorem~\ref{the:extcomp}
we may now assume that we are not in the case
discussed in the previous
subsection.
Hence we may now assume that we
have exactly two components of maximal
width~$m$ and two components of maximal
length~$\tilde m=|\mu\cap\nu'|$;
we know that these are four different components
as no constituent is both of maximal
width and length by Theorem~\ref{PConst}.
Our task is to find a fifth component in the product $[\mu]\cdot[\nu]$
which is almost extreme.
}

We set $\tga=\mu\cap \nu'$.
Since we have two components of maximal width and length, respectively,
both products $[\mu / \ga ]\cdot[\nu / \ga]$
and $[\mu / \tga ]\cdot[\nu / \tga]$ have exactly two components.
{Note that here $d>1$.}

We focus on the first product.
This situation splits into the following cases.

\begin{enumerate}
\item $[\mu / \gamma ]=[1^2]+[2]$ and $[\nu / \gamma]=[1^2]+[2]$.
\item $[\mu / \ga ]$ and $[\nu / \ga]$ are both irreducible, and the product has two components.
\item $[\mu / \ga ]$ has two components and $[\nu / \ga]$ is of degree~1.
\end{enumerate}

\bigskip

{\bf 6.2.1.} First we treat Case~(1),
where both skew diagrams decompose into two disconnected nodes.
Note that the assumptions imply that $n\geq 6$ and that $\ga$ has at least
three removable nodes.

By Corollary~\ref{CDvir2}(i) we obtain in this case from
$[\mu / \ga ]\cdot [\nu / \ga ] = 2[2]+2[1^2]$
the constituents  $[n-2,2]$ and $[n-2,1^2]$ of maximal width
in $[\mu]\cdot [\nu]$.

We will show that all three possible components of almost maximal width~$n-3$
appear in the product, using Lemma~\ref{lem:prepalmext}.
Assume there is a removable $\ga$-node $A_0$ which
is disconnected from $\mu/\ga$; then
$[\mu / \ga_{A_0} ]=[3]+2[2,1]+[1^3]$.
We want to compute
$$
\chi_0=[\mu / \ga_{A_0} ]\cdot [\nu / \ga_{A_0} ]
- ([\mu / \ga ]\cdot [\nu / \ga ])\uparrow
=
[\mu / \ga_{A_0} ]\cdot [\nu / \ga_{A_0} ] - (2[3]+4[2,1]+2[1^3])\:.
$$

If $A_0$ is connected to at most one of the nodes of $\nu/\ga$, then
$[\nu/\ga_{A_0}]$ contains $[3]+[2,1]$ or $[2,1]+[1^3]$,
and thus
$$
\chi_0'= (3[3]+6[2,1]+3[1^3]) - (2[3]+4[2,1]+2[1^3])
= [3]+2[2,1]+[1^3]
$$
is a character contained in the character~$\chi_0$.
Since $n\geq 6$, we get three constituents
$[n-3,3]$, $[n-3,2,1]$, $[n-3,1^3]$ in the product $[\mu]\cdot[\nu]$.

The same argument can be used with $\mu,\nu$ interchanged.
Hence we may now assume that every removable $\ga$-node
is connected to two of the four nodes of $\mu/\ga$ and~$\nu/\ga$.
Arguing from the top removable $\ga$-node down, one easily sees
that then we must have $\ga=(3,2,1)$.
If every removable $\ga$-node is connected to both a $\mu/\ga$ and a $\nu/\ga$-node,
then
$$
\begin{array}{rcl}
\chi &=&
\dstyle \sum_{A \ \ga\text{-node}} [\mu / \ga_{A} ]\cdot [\nu / \ga_{A} ]
- ([\mu / \ga ]\cdot [\nu / \ga ])\uparrow \\[8pt]
&=& 3([2,1]+[3])\cdot ([2,1]+[1^3]) - (2[3]+4[2,1]+2[1^3])
\\[8pt]
&=& [3]+5[2,1]+4[1^3]\:.
\end{array}
$$
If some removable $\ga$-node $A_0$ is connected to two nodes of the same skew diagram
$\mu/\ga$ or~$\nu/\ga$,
then there is also at least one removable $\ga$-node $A_1$ that is
connected to both $\mu/\ga$ and~$\nu/\ga$,
and we have
$$
\begin{array}{rcl}
\chi'&=&
\dstyle \sum_{A \in \{A_0,A_1\}} [\mu / \ga_{A} ]\cdot [\nu / \ga_{A} ]
- ([\mu / \ga ]\cdot [\nu / \ga ])\uparrow \\[8pt]
&=& [2,1]\cdot ([3]+2[2,1]+[1^3]) + ([2,1]+[3])\cdot ([2,1]+[1^3]) - (2[3]+4[2,1]+2[1^3])
\\[8pt]
&=& [3]+3[2,1]+2[1^3] \:,
\end{array}
$$
a subcharacter of the character~$\chi$.
Hence in both cases we get all three almost extreme constituents
$[n-3,3]$, $[n-3,2,1]$, $[n-3,1^3]$ in the product $[\mu]\cdot[\nu]$.

{
Thus in any case we have found at least five (almost) extreme components in the product.}

\bigskip

{\bf 6.2.2.} We now deal with Case (2).
By  Theorem~\ref{Th2C} we know that up to conjugation and renaming
we are in the following situation:\\
\centerline{$[\mu / \ga ]=[\al]$, $[\nu / \ga]=[d-1,1]$, where
$\al=(a^b)$ is a nontrivial rectangle, i.e., $a,b>1$.}

By Proposition~\ref{prop:5.5}
we then obtain at least five components of almost maximal width $m-1$
in $[\mu]\cdot [\nu]$.

\bigskip

{\bf 6.2.3.}
Now to Case (3), where
$[\mu / \ga ]$ has two components and $[\nu / \gamma]$ is of degree 1.
We may assume that $\nu/\ga$ is a row.

By Proposition~\ref{SkCh2Comp}, we know that the
skew character  $[\mu/\ga]$ has the form
$$ [\mu/\ga] =[\al^X]+[\al^Y]$$
for some partition $\al$ and two distinct addable nodes $X$, $Y$ for~$\al$.

Thus $[\mu]\cdot [\nu]$ has $[m,\al^X]$ and $[m,\al^Y]$
as its components of maximal width.

Using Lemma~\ref{ThDvir2'}, we want to produce components of almost maximal
width $m-1$ in the product.
Thus we now consider
$$
\chi=\sum_{A\ \ga\mathrm{-node}}[\mu / \gamma_A] \cdot [\nu / \gamma_A]
- [\mu/ \ga]\uparrow \:.$$

First we assume that there is a removable $\ga$-node $A_0$ disconnected from $\mu/\ga$.
Then
$$
[\mu/ \ga_{A_0}]=[\mu/ \ga]\uparrow = [\al^X]\uparrow+[\al^Y]\uparrow \:.
$$
If $[d,1]$ is a constituent of $[\nu/\ga_{A_0}]$, then
$
\chi_0=[\mu / \ga_{A_0}] \cdot [\nu / \ga_{A_0}]
- [\mu/ \ga]\uparrow $ contains
$$
\chi_0'=([\mu / \ga]\uparrow ) \cdot [d,1]
- [\mu/ \ga]\uparrow
= ([\mu / \ga] \cdot [d-1,1])\uparrow \:.
$$
Now $([\al^X]+[\al^Y]) \cdot [d-1,1]$ clearly
contains $[\al^Y]+[\al^X]=[\mu/\ga]$ as a subcharacter, hence
the character $\chi_0$ contains
$$\chi_0''=[\mu/\ga]\uparrow = [\al^X]\uparrow+[\al^Y]\uparrow \:.
$$
If one of $\al^X$ and $\al^Y$ is not a rectangle,
we  clearly have $c(\chi) \geq c(\chi_0'') \geq 3$.
That $\al^X$ and $\al^Y$ are both rectangles can only occur when
$\al=(1)$, i.e., when $[\mu/\ga]=[2]+[1^2]$.
In this case we have
$$
\chi_0''= ([2]+[1^2])\uparrow = [3]+2[2,1]+[1^3]\:,
$$
hence again  $c(\chi) \geq 3$.
Thus in any case
we have three components of width $m-1$ in~$[\mu]\cdot[\nu]$.

We may now assume that $\nu/\ga_{A_0}$ is a row, and furthermore,
that any removable $\ga$-node $A_1\neq A_0$ is connected to $\mu/\ga$.
Then
$$
\chi_0=[\mu / \ga_{A_0}] \cdot [\nu / \ga_{A_0}]
- [\mu/ \ga]\uparrow =0\:.$$
Since $\mu/\ga$ is a proper skew diagram, there must be
such a $\ga$-node $A_1$ which is an inner node for
$\mu/\ga$, i.e., it is connected to a node of $\mu/\ga$
but it is not above the highest row nor to the left of
the leftmost column of $\mu/\ga$.
In fact, considering the list in Proposition~\ref{SkCh2Comp}
we see that then
$$
[\mu/\ga_{A_1}]= \left\{
\begin{array}{ll}
[\al^{XY}] & \text{if $\mu/\ga_{A_1}$ is a partition}\\
\-[\al^{XX'}] + [\al^{YY'}] & \text{if $\mu/\ga_{A_1}$ is not a partition but
on the list in~Prop.~\ref{SkCh2Comp} }
\end{array}
\right.
$$
where $X'$, $Y'$ are addable nodes for $\al^X$, $\al^Y$, respectively
(possibly $X'=Y$ or $Y'=X$, but not both);
note that $X$, $Y$ are nodes that we can add
in an independent way, whereas $X$, $X'$ and $Y$, $Y'$ may only
be added in this order. When we are not in one of the two
situations above, $[\mu/\ga_{A_1}]$ is a skew character
with at least three components, including $[\al^{XX'}] + [\al^{YY'}]$
as above.

Now whenever $A_1$ is disconnected from $\nu/\ga$ we have
$$
\chi_1=[\mu/\ga_{A_1}]\cdot [\nu/\ga_{A_1}]=
[\mu/\ga_{A_1}]\cdot ([d]\uparrow) =
[\mu/\ga_{A_1}]\downarrow\uparrow \:.
$$
>From the description above, we see immediately that
$[\mu/\ga_{A_1}]\downarrow$
contains $[\al^X]+[\al^Y]=[\mu/\ga]$, hence arguing as above,
$\chi_1$ (and thus $\chi$)
has at least three components.

Thus we assume now that all removable $\ga$-nodes that are inner nodes
for $\mu/\ga$ are connected to $\nu/\ga$; this can only happen
if $\mu/\ga$ is a disconnected skew diagram with two parts and
$\nu/\ga$ between them.
There can only be one such $\ga$-node $A_1\neq A_0$.
We then have
$$
\chi_1=[\mu/\ga_{A_1}]\cdot [\nu/\ga_{A_1}]=
[\mu/\ga_{A_1}]\cdot [d,1] \:.
$$
>From the description above, we deduce that $\chi_1$ always has
$[\al^{XY}]$ as a constituent, giving a component
$[m-1,\al^{XY}]$ in $[\mu]\cdot[\nu]$ which is
of almost maximal width but not of maximal length
as $\ell(\al^{XY})=\max(\ell(\al^X),\ell(\al^Y))$
and there is no constituent which is maximal in both respects.

\smallskip

We may now assume that all removable $\ga$-nodes
are connected to $\mu/\ga$.

Let  $A_0$ be a removable $\ga$-node, and
assume that this is not connected to $\nu/\ga$.
In any case, $[\mu/\ga_{A_0}]\downarrow$
contains $[\al^X]+[\al^Y]=[\mu/\ga]$, so that
$$\chi_0
= [\mu/\ga_{A_0}]\cdot [\nu/\ga_{A_0}]  - [\mu/\ga]\uparrow
= [\mu/\ga_{A_0}]\downarrow\uparrow  - [\mu/\ga]\uparrow
$$
is a character (or 0).

Since $\mu/\ga$ is a proper skew character,
there has to be a second removable $\ga$-node $A_1\neq A_0$.
If $[\nu/\ga_{A_1}]$ contains $[d,1]$, then as before,
$\chi_1
= [\mu/\ga_{A_1}]\cdot [\nu/\ga_{A_1}]$
contains $[\al^{XY}]$, and this provides a fifth
almost extreme component $[m-1,\al^{XY}]$ in $[\mu]\cdot[\nu]$
as we have seen above.
Thus we may now assume that $[\nu/\ga_{A_1}]=[d+1]$.

Now if $A_1$ is an inner node for $\mu/\ga$, then
$[\mu/\ga_{A_1}]$ always has a constituent
of length $\max(\ell(\al^X),\ell(\al^Y))$, namely
the one coming from sorting the rows, and this
would also provide a (fifth) almost extreme component.

Thus $A_1$ can only be an outer node for $\mu/\ga$.
Proposition~\ref{SkCh2Comp} implies that
$[\mu/\ga_{A_1}]$ has at least three components,
except in the case where $[\mu/\ga]=[1^{a+1}]\otimes [r]$
and we obtain
$[\mu/\ga_{A_1}]=[1^{b+1}]\otimes [s]$ with
either $b=a+1$ or $s=r+1$.
But in this case, $[\mu/\ga_{A_1}]$ has
$[\al^{XY}]$ as a constituent, giving a (fifth) almost extreme
component  $[m-1,\al^{XY}]$ in $[\mu]\cdot[\nu]$ as before.

Hence we may now assume that any removable $\ga$-node
is connected to  $\nu/\ga$.
Since $\mu/\ga$ is a proper skew diagram and we always
have a removable $\ga$-node that is inner with respect to
$\mu/\ga$, this can only be true when $\mu/\ga$ is disconnected,
i.e., we are in one of the cases (i) or (ii) of Proposition~\ref{SkCh2Comp}
(as before, up to translation
and order of the two connected parts).
Furthermore, $\nu/\ga$ sits between the two
parts of $\mu/\ga$ and we must have two
removable $\ga$-nodes $A_0$, $A_1$ connected to both
$\mu/\ga$ and~$\nu/\ga$.
Altogether, there are now only four cases we have to consider.

Here are the corresponding pictures (we mark both parts of $\mu/\ga$ by $\mu/\ga$):
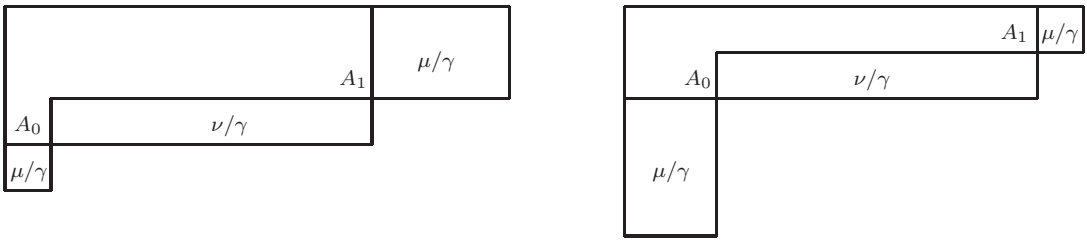
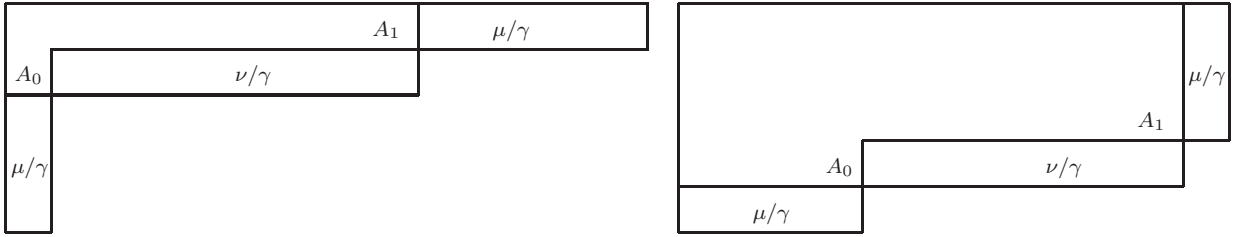
\begin{center}
\setlength{\unitlength}{0.024in}
\begin{picture}(0,30)(110,50)
\thicklines
\put( 0,80){\line( 1, 0){110}}
\put( 10,60){\line( 1, 0){100}}
\put( 0,50){\line( 1, 0){80}}
\put( 0,40){\line( 1, 0){ 10}}
\put( 0,80){\line( 0,-1){ 40}}
\put(10,60){\line( 0,-1){ 20}}
\put(80,80){\line( 0,-1) {30}}
\put(110,80){\line( 0,-1){ 20}}
\put(90,67){\makebox(0,0)[lb]{\smash{\SetFigFont{8}{7.2}{\rmdefault
}{\mddefault}{\updefault}$\mu /\gamma$}}}
\put(-5,43){
\makebox(0,0)[lb]{
\smash{\SetFigFont{8}{7.2}{\rmdefault}{\mddefault}{\updefault}
$\mu /\gamma$}
}}
\put(45,53){\makebox(0,0)[lb]{\smash{\SetFigFont{8}{7.2}{\rmdefault}{\mddefault}{\updefault}$\nu /\gamma$}}}
\put(73,63){\makebox(0,0)[lb]{\smash{\SetFigFont{8}{7.2}{\rmdefault}{\mddefault}{\updefault}$A_1$}}}
\put(2,53){\makebox(0,0)[lb]{\smash{\SetFigFont{8}{7.2}{\rmdefault}{\mddefault}{\updefault}$A_0$}}}
\end{picture}
\setlength{\unitlength}{0.024in}%
\hfil
\begin{picture}(-100,50)(100,50)
\thicklines
\put( 0,80){\line( 1, 0){100}}
\put( 20,70){\line( 1, 0){80}}
\put( 0,60){\line( 1, 0){90}}
\put( 0,30){\line( 1, 0){ 20}}
\put( 0,80){\line( 0,-1){ 50}}
\put(20,70){\line( 0,-1){ 40}}
\put(90,80){\line( 0,-1) {20}}
\put(100,80){\line( 0,-1){ 10}}
\put(87,73){\makebox(0,0)[lb]{
\smash{\SetFigFont{8}{7.2}{\rmdefault}{\mddefault}{\updefault}
$\mu /\gamma$}
}}
\put(0,43){
\makebox(0,0)[lb]{
\smash{\SetFigFont{8}{7.2}{\rmdefault}{\mddefault}{\updefault}
$\mu /\gamma$}
}}
\put(50,63){\makebox(0,0)[lb]{\smash{\SetFigFont{8}{7.2}{\rmdefault}{\mddefault}{\updefault}$\nu /\gamma$}}}
\put(82,73){\makebox(0,0)[lb]{\smash{\SetFigFont{8}{7.2}{\rmdefault}{\mddefault}{\updefault}$A_1$}}}
\put(13,63){\makebox(0,0)[lb]{\smash{\SetFigFont{8}{7.2}{\rmdefault}{\mddefault}{\updefault}$A_0$}}}
\end{picture}

\setlength{\unitlength}{0.024in}%
\begin{picture}(10,110)(115,10)
\thicklines
\put( 0,80){\line( 1, 0){140}}
\put( 10,70){\line( 1, 0){130}}
\put( 0,60){\line( 1, 0){90}}
\put( 0,30){\line( 1, 0){ 10}}
\put( 0,80){\line( 0,-1){ 50}}
\put(10,70){\line( 0,-1){ 40}}
\put(90,80){\line( 0,-1) {20}}
\put(140,80){\line( 0,-1){ 10}}
\put(100,73){
\makebox(0,0)[lb]{
\smash{\SetFigFont{8}{7.2}{\rmdefault}{\mddefault}{\updefault}
$\mu /\gamma$}
}}
\put(-5,43){
\makebox(0,0)[lb]{
\smash{\SetFigFont{8}{7.2}{\rmdefault}{\mddefault}{\updefault}
$\mu /\gamma$}
}}
\put(50,63){\makebox(0,0)[lb]{\smash{\SetFigFont{8}{7.2}{\rmdefault}{\mddefault}{\updefault}$\nu /\gamma$}}}
\put(80,73){\makebox(0,0)[lb]{\smash{\SetFigFont{8}{7.2}{\rmdefault}{\mddefault}{\updefault}$A_1$}}}
\put(2,63){\makebox(0,0)[lb]{\smash{\SetFigFont{8}{7.2}{\rmdefault}{\mddefault}{\updefault}$A_0$}}}
\end{picture}
\hfil
\setlength{\unitlength}{0.024in}%
\begin{picture}(-100,80)(100,10)
\thicklines
\put( 0,80){\line( 1, 0){120}}
\put( 40,50){\line( 1, 0){80}}
\put( 0,40){\line( 1, 0){110}}
\put( 0,30){\line( 1, 0){ 40}}
\put( 0,80){\line( 0,-1){ 50}}
\put(40,50){\line( 0,-1){ 20}}
\put(110,80){\line( 0,-1) {40}}
\put(120,80){\line( 0,-1){ 30}}
\put(105,63){
\makebox(0,0)[lb]{
\smash{\SetFigFont{8}{7.2}{\rmdefault}{\mddefault}{\updefault}
$\mu /\gamma$}
}}
\put(10,33){
\makebox(0,0)[lb]{
\smash{\SetFigFont{8}{7.2}{\rmdefault}{\mddefault}{\updefault}
$\mu /\gamma$}
}}
\put(80,43){\makebox(0,0)[lb]{\smash{\SetFigFont{8}{7.2}{\rmdefault}{\mddefault}{\updefault}$\nu /\gamma$}}}
\put(100,53){\makebox(0,0)[lb]{\smash{\SetFigFont{8}{7.2}{\rmdefault}{\mddefault}{\updefault}$A_1$}}}
\put(32,43){\makebox(0,0)[lb]{\smash{\SetFigFont{8}{7.2}{\rmdefault}{\mddefault}{\updefault}$A_0$}}}
\end{picture}
\end{center}

In all cases, we let $A_0$ and $A_1$ be the removable $\ga$-nodes
such that $\nu/\ga_{A_0}$ is a row
and $[\nu/\ga_{A_1}]=[d,1]$.
We know that  $[\mu/\ga]=[\al^X]+[\al^Y]$,
with  $X$ being the top and $Y$ the bottom
addable node in our situation,
where now $\al$ is a nontrivial rectangle
or a hook.
Let $X'$ be the top addable node of $\al^X$.
One easily checks that then
in all cases, $[\mu/\ga_{A_0}]$ has a constituent
$[\al^{XY}]$ and
 $[\mu/\ga_{A_1}]$ has at least the constituents
$[\al^{XY}]$ and $[\al^{XX'}]$.
Hence, similarly as before,
$[\mu/\ga_{A_1}]\cdot [d,1]$ contains
$[\al^Y]\uparrow+[\al^X]\uparrow=[\mu/\ga]\uparrow$.
Thus in all cases, we have the following:
$$\chi =  [\mu/\ga_{A_0}] + [\mu/\ga_{A_1}]\cdot [d,1] - [\mu/\ga]\uparrow$$
is a character containing $[\al^{XY}]$,
and this produces, as before, a constituent $[m-1,  \al^{XY}]$
in $[\mu]\cdot[\nu]$ which is of almost maximal width and
not of maximal length.

{
Hence we have found in all cases of a product with exactly
two components of maximal width and two of maximal length
a further almost extreme component, and thus we are done.
}
\qed

\bigskip

Finally we turn to the confirmation of the classification conjecture
in~\cite{BK}.
We recall the classification result we want to prove below;
of course, the decompositions of
all the products appearing in the statements
below are known.

\begin{theorem}
Let $\mu$, $\nu \vdash n$.
\begin{enumerate}
\item[(i)] \label{Th3C}\label{the:main3}
We have $c([\mu]\cdot [\nu])=3$
if and only if $n=3$ and $\mu=\nu=(2,1)$ or $n=4$ and $\mu=\nu=(2,2)$.

The product is then one of
\\
\centerline{$[2,1]^2 = [3] + [2,1] + [1^3]$\:,
$[2^2]^2 = [4] + [2^2] + [1^4]$\:.}
\item[(ii)]\label{Th4C}\label{the:main4}
We have $c([\mu]\cdot [\nu])=4$
if and only if one of the following holds:
\begin{enumerate}
\item[(1)] $n\geq 4$ and $\mu,\nu\in\{(n-1,1),(2,1^{n-2})\}$;
here the products are
$$\begin{array}{rcl}
[n-1,1]^2 =[2,1^{n-2}]^2 &=& [n] + [n-1,1]+ [n-2,2] + [n-2,1^2]\\[5pt]
\-
[n-1,1]\cdot[2,1^{n-2}] &=& [1^n] + [2,1^{n-2}]+ [2^2,1^{n-4}] + [3,1^{n-3}]\
\, .
\end{array}$$
\item[(2)] $n=2k+1$ for some $k\geq 2$, and
one of $\mu$, $\nu$ is in $\{(2k,1),(2,1^{2k-1})\}$
while the other one is in $\{(k+1,k),(2^k,1)\}$;
here the products are
$$
\begin{array}{rcl}
[2k,1]\cdot [k+1,k] &=&
[k+2,k-1] + [k+1,k] + [k+1,k-1,1] + [k^2,1]\\[5pt]
\-
[2,1^{2k-1}] \cdot [k+1,k] &=&
[2^{k-1},1^3] + [2^k,1] + [3,2^{k-2},1^2] + [3,2^{k-1}]\, .
\end{array}
$$
\item[(3)] $n=6$ and $\mu,\nu\in\{(2^3),(3^2)\}$;
here we have
$$[3^2]^2=[6]+[4,2]+[3,1^3]+[2^3]\,, \:
[3^2]\cdot[2^3]=[1^6]+[2^2,1^2]+[4,1^2]+[3^2]\:.$$
\end{enumerate}
\end{enumerate}
\end{theorem}
\medskip

\proof
We only need to check the exceptional cases where
a product has at most four (almost) extreme components
listed in Theorem~\ref{the:extcomp}.

In Section~\ref{sec:specprods}, we have already classified all
products with the character $[n-1,1]$ that have at most four
components in
Corollary~\ref{Cor:Nat},
while the classification of squares $[\la]^2$ with at most four components
is stated in Proposition~\ref{LTensSq}.
The exceptional case of the product $[3^2]\cdot [4,2]$
appeared in the investigation of 2-part partitions in Proposition~\ref{twopart},
where we already noticed that it has five components.

This information allows to handle all cases in Theorem~\ref{the:extcomp},
and we get exactly the classification in Theorem~\ref{the:34C}.
\qed


\end{document}